\documentclass[a4paper,11pt]{article}

\usepackage[top=2.5cm,bottom=2cm,left=2.6cm,right=2.6cm]{geometry}

\usepackage[english]{babel}
\usepackage[utf8x]{inputenc}
\usepackage[T1]{fontenc}
\usepackage{bold-extra}

\usepackage{xcolor}
\usepackage{tikz}
\usetikzlibrary{arrows.meta}
\usepackage{graphicx}
\usepackage{subcaption}
\usepackage{multirow, multicol}
\usepackage{authblk}
\usepackage{amsmath,amsfonts,amssymb,amsthm}
\usepackage{natbib}
\usepackage[colorlinks=true, allcolors=blue]{hyperref}
\usepackage{xspace}

\usepackage[abs]{overpic}% option [abs] for absolute values 

\usepackage{lipsum}

\usepackage{algcompatible}
\usepackage{algorithm}
\usepackage[noend]{algpseudocode}
\floatname{algorithm}{Procedure}

\renewcommand{\algorithmiccomment}[1]{// #1}
\algnewcommand{\LeftComment}[1]{\Statex \algorithmiccomment{#1}}
\algnewcommand{\AlgSection}[1]{\item[] \Statex ==== \textbf{#1} ====}
\algnewcommand{\AlgSubsection}[1]{\item[] \Statex --- \textbf{#1} ---}
\makeatletter
\def\BState{\State\hskip-\ALG@thistlm}
\makeatother

\usepackage{mathtools}
\newcommand{\Set}[1]{\mathcal{#1}}
\newcommand{\Param}[1]{#1}

\providecommand{\card}[1]{\lvert#1\rvert}

\newcommand{\fac}{\mathcal{F}}
\newcommand{\ware}{\mathcal{W}}
\newcommand{\outl}{\mathcal{O}}
\newcommand{\sku}{\mathcal{S}}
\newcommand{\pack}{\mathcal{P}}
\newcommand{\mov}{\mathcal{M}}
\newcommand{\sini}{IS}
\newcommand{\fixd}{FD}
\newcommand{\vard}{VD}
\newcommand{\pri}{\pi}
\newcommand{\wei}{W}
\newcommand{\capa}{Cap}
\newcommand{\cost}{Cost}
\newcommand{\yvar}{\boldsymbol{Y}}
\newcommand{\xvar}{\boldsymbol{X}}
\newcommand{\ypack}{\boldsymbol{y}}
\newcommand{\xpack}{\boldsymbol{x}}
\newcommand{\sfin}{\boldsymbol{FS}}

\newcommand{\avap}{AP}

\newcommand{\problem}[1]{{\sc \textbf{#1}$^\text{{\sc{\textbf{p}}}}$}\xspace}
\newcommand{\problemm}[1]{{\sc \textbf{#1}$^\text{{\sc{\textbf{p}}}}$}}

\newcommand{\tp}{\problem{t}}
\newcommand{\rtp}{{\sc \textbf{rt}$_\delta^\text{{\sc{\textbf{p}}}}$}\xspace}

\newcommand{\rp}{\problem{r}}
\newcommand{\pp}{\problem{p}}
\newcommand{\ppm}[1]{{\sc \textbf{p}$_{#1}^\text{{\sc{\textbf{p}}}}$}}
\newcommand{\rpm}[1]{{\sc \textbf{r}$_{#1}^\text{{\sc{\textbf{p}}}}$}}
\newcommand{\tpnosp}{\problemm{t}}
\newcommand{\rtpnosp}{{\sc \textbf{rt}$_\delta^\text{{\sc{\textbf{p}}}}$}}
\newcommand{\rpnosp}{\problemm{r}}

\newcommand{\tppp}{\tpnosp\texttt{+}\pp}
\newcommand{\rtprp}{\rtpnosp\texttt{+}\rp}
\newcommand{\rtprppp}{\rtpnosp\texttt{+}\rpnosp\texttt{+}\pp}
\newcommand{\TRPA}{\solver{RT-R-P}\xspace}
\newcommand{\TRPAdelta}[1]{\solver{RT$_{#1}$-R-P}\xspace}
\newcommand{\TPA}{\solver{T-P}\xspace}

\usepackage{enumitem}
\newcounter{myitemcounter}

\newcommand{\myitemlabel}{$\bullet$\ }

\newcommand{\myitem}{%
\stepcounter{myitemcounter}
\myitemlabel
}

\newcommand{\anitem}[1]{%
\myitem #1 &
}

\newcommand{\lastitem}[1]{%
\myitem #1 \\
}

\newenvironment{inlineitemize}
{\setcounter{myitemcounter}{0}
\begin{tabular}{llllllllll} % you won't want more columns
}
{\end{tabular}}

\title{Optimal policies for stock redistribution in a retail network: Mathematical modeling and algorithmic solution\thanks{\textbf{Acknowledgments.} The authors are grateful to Nanos' staff for fruitful discussions on the specifics of their daily operations. This work is part of the R\&D projects PID2021-124030NB-C31 and PID2021-124030NB-C32 funded by MICIU/AEI/10.13039/501100011033/ and by ERDF/EU. This research was also funded by Grupos de Referencia Competitiva ED431C-2021/24 from the Consellería de Cultura, Educación e Universidades, Xunta de Galicia. We wish to acknowledge the support received from the Centro de Investigaci\'on de Galicia "CITIC", funded by Xunta de Galicia and the European Union (European Regional Development Fund- Galicia 2014-2020 Program), by grant ED431G 2019/01.}}
\author[1,2]{Julio Gonz\'alez-D\'iaz}
\author[1,2]{\'Angel M. Gonz\'alez-Rueda\thanks{Corresponding author: angelmanuel.gonzalez.rueda@usc.es.}}
\author[1]{Irene Llana Garc\'ia}
\author[1]{Jorge Rodr\'iguez Veiga}
\affil[1]{Department of Statistics, Mathematical Analysis and Optimization and MODESTYA Research Group. University of Santiago de Compostela. Santiago de Compostela, Spain.}
\affil[2]{CITMAga (Galician Center for Mathematical Research and Technology), Santiago de Compostela, Spain.}

\newcommand{\solver}[1]{\texttt{#1}}

\date{\today}

\begin{document}

\maketitle

\begin{abstract}
We study the problem of stock replenishment and transshipment in the retail industry. We develop a model that can accommodate different policies, including centralized redistribution (replenishment) and decentralized redistribution (lateral transshipments), allowing for direct comparisons between them. We present a numeric analysis in which the benchmark instances stem from the collaboration with a high-end clothing retail company. The underlying model, as usually in the field, is a large-scale mixed-integer linear programming problem. We develop a specific algorithmic procedure to solve this MILP problem and compare its performance with the direct solution via state-of-the-art solvers. 
\end{abstract}

\textbf{Keywords} Suppy chain, inventory management, stock redistribution, lateral transshipments, MILP problems, retail industry.

\section{Introduction}
Efficient inventory management is at the core of most modern companies. In this paper we study the specific issue of stock redistribution in a two-echelon supply chain network, with warehouses and outlets in the first and second echelons, respectively. The models we develop have been driven by a collaboration with Nanos, a company that manufactures and sells high-end kid's clothing.\footnote{See \url{https://www.nanos.es/us/}.} In general, the most common reason to need some form of redistribution is a stockout, which is often met with a backorder demand to a warehouse or distribution center. In the two-echelon setting under study, since there is a network of facilities, one can benefit from inventory pooling via (lateral) transshipments from other outlets. As discussed in \cite{Paterson2011}, ``\textit{These transshipments can be conducted periodically at predetermined points in time to proactively redistribute stock, or they can be used reactively as a method of meeting demand which cannot be satisfied from stock on hand}''. Thus, redistribution does not only arise as a response to the so called ``emergency transshipments'' (reactive, as in the pioneering work by \cite{Krishnan1965}), but may be a regular operation that is performed several times within a selling season to balance the stocks at the different outlets given past sales and estimates for future demand (proactive, as in the seminal papers by \cite{Allen1958,Allen1961,Allen1962} and \cite{Gross1963}).\footnote{In the particular case of the clothing industry, redistribution can also be used to consolidate the available sizes of a given article at some of the outlets, as discussed in \cite{Caro2010b}.} Although the mathematical model we develop is motivated by the proactive redistribution case, it can also be applied to the reactive one and to combinations of both. 

In our approach we take as given a snapshot of the whole network on which an optimal redistribution is sought. This snapshot is the input of our optimization model and includes, for each stock keeping unit (SKU),\footnote{In the clothing industry, given a specific garment, the SKU has to account for all attributes that distinguish it from other types of the same garment, such as size, color, and material.} the available stock at the warehouses and the available stock along with estimates for future demand at each outlet. Transportation costs, which may differ across pairs of facilities, are also known. This paper has two main contributions to the efficient management of networks in which the aforementioned information is available: a novel redistribution model and a tailor-made algorithmic procedure to solve it.

The first contribution consists of the development in Section~\ref{sec:model} of a model to compute the optimal redistribution strategy, which should provide the right balance between replenishments from warehouses and lateral transshipments. This model allows to evaluate the relative efficiency of special policies in which some forms of replenishment or transshipment are not available or just	 partially limited. More generally, the model can also be used to assess the convenience of potential network expansions (new outlets, warehouses or distribution centers) given their impact, for instance, on the usual redistribution scenarios.

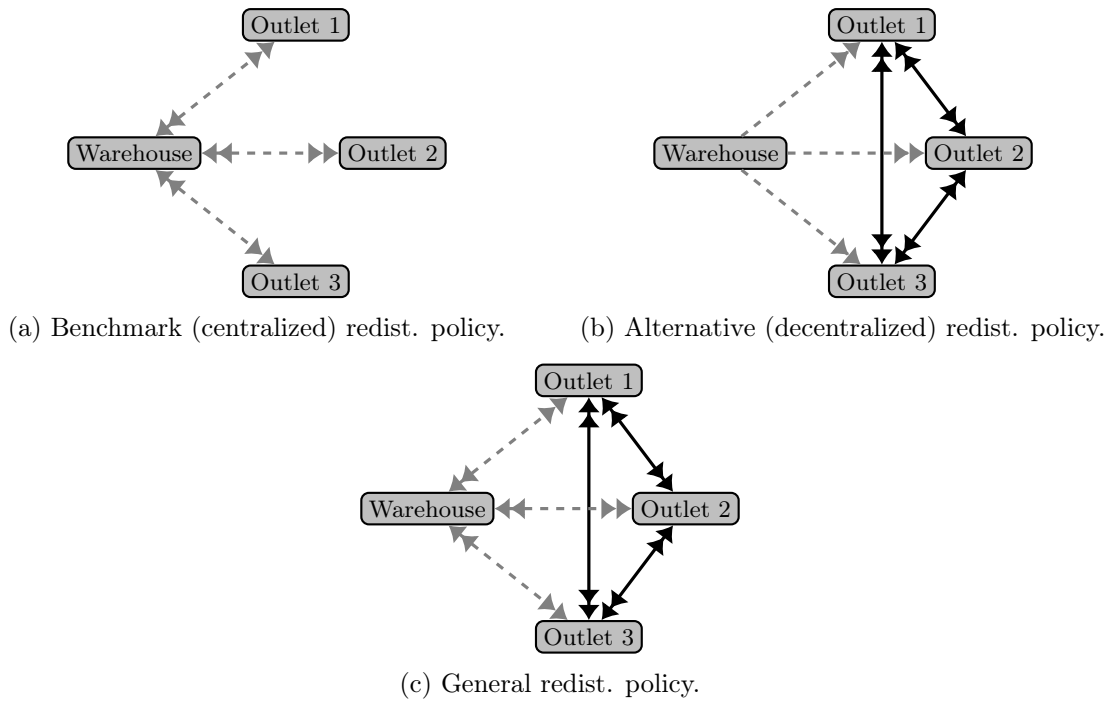
\begin{figure}[!htbp]
\centering
\begin{subfigure}{.48\textwidth}
  \centering
	\begin{tikzpicture}[scale=0.85,>={Latex[length=2mm,width=3mm]},font=\footnotesize,every node/.style={font=\footnotesize},thick,estilo/.style={rectangle,draw,inner sep=1mm,rounded corners=1mm,fill=gray!50},theshift/.style={xshift=7.5cm,yshift=0cm}]
		\node[estilo] (w1) at (-1.5,0){Warehouse};
		\node[estilo] (o1) at (1,2){Outlet 1};
		\node[estilo] (o2) at (2.5,0){Outlet 2};
		\node[estilo] (o3) at (1,-2){Outlet 3};
		%\node[estilo] (o4) at (1.5,-3){Outlet 4};
		
		%\draw[<<->>,very thick] (o1) -- (o2);
		%\draw[<<->>,very thick] (o2) -- (o3);
		%\draw[<<->>,very thick] (o3) -- (o4);
		%\draw[<<->>,very thick] (o4) to[bend right=90] (o1);
		%\draw[<<->>,very thick] (o4) to[bend right=60] (o2);
		%\draw[<<->>,very thick] (o3) -- (o1);
		
		\draw[<<->>,very thick, dashed, gray] (w1) -- (o1);
		\draw[<<->>,very thick, dashed, gray] (w1) -- (o2);
		\draw[<<->>,very thick, dashed, gray] (w1) -- (o3);
		%\draw[<<->>,very thick, dashed, gray] (w1) -- (o4);
		
	\end{tikzpicture}
	\caption{Benchmark (centralized) redist. policy.}
  \label{subfig:centralized}
\end{subfigure}
\vspace{0.2cm}
\begin{subfigure}{.48\textwidth}
  \centering
	\begin{tikzpicture}[scale=0.85,>={Latex[length=2mm,width=3mm]},font=\footnotesize,every node/.style={font=\footnotesize},thick,estilo/.style={rectangle,draw,inner sep=1mm,rounded corners=1mm,fill=gray!50},theshift/.style={xshift=7.5cm,yshift=0cm}]
		\node[estilo] (w1) at (-1.5,0){Warehouse};
		\node[estilo] (o1) at (1,2){Outlet 1};
		\node[estilo] (o2) at (2.5,0){Outlet 2};
		\node[estilo] (o3) at (1,-2){Outlet 3};
		
		\draw[<<->>,very thick] (o1) -- (o2);
		\draw[<<->>,very thick] (o2) -- (o3);
		\draw[<<->>,very thick] (o3) -- (o1);
		
		\draw[->>,very thick, dashed, gray] (w1) -- (o1);
		\draw[->>,very thick, dashed, gray] (w1) -- (o2);
		\draw[->>,very thick, dashed, gray] (w1) -- (o3);
	\end{tikzpicture}
	\caption{Alternative (decentralized) redist. policy.}
  \label{subfig:decentralized}
\end{subfigure}
\hspace{0.1cm}
\begin{subfigure}{.48\textwidth}
  \centering
	\begin{tikzpicture}[scale=0.85,>={Latex[length=2mm,width=3mm]},font=\footnotesize,every node/.style={font=\footnotesize},thick,estilo/.style={rectangle,draw,inner sep=1mm,rounded corners=1mm,fill=gray!50},theshift/.style={xshift=7.5cm,yshift=0cm}]
		\node[estilo] (w1) at (-1.5,0){Warehouse};
		\node[estilo] (o1) at (1,2){Outlet 1};
		\node[estilo] (o2) at (2.5,0){Outlet 2};
		\node[estilo] (o3) at (1,-2){Outlet 3};
		
		\draw[<<->>,very thick] (o1) -- (o2);
		\draw[<<->>,very thick] (o2) -- (o3);
		\draw[<<->>,very thick] (o3) -- (o1);
		
		\draw[<<->>,very thick, dashed, gray] (w1) -- (o1);
		\draw[<<->>,very thick, dashed, gray] (w1) -- (o2);
		\draw[<<->>,very thick, dashed, gray] (w1) -- (o3);
	\end{tikzpicture}
	\caption{General redist. policy.}
  \label{subfig:general}
\end{subfigure}

	\caption{Different redistribution policies captured by the model developed in Section~\ref{sec:model}.}
	\label{fig:redist_schemes}
\end{figure}

In Figure~\ref{fig:redist_schemes} we represent three different policies captured by our model:\footnote{For the sake of illustration, Figure~\ref{fig:redist_schemes} represents just one central warehouse, but our model allows for multiple warehouses.}
\begin{itemize}
	\item Figure~\ref{subfig:centralized} represents the redistribution policy in place at Nanos at the time of this study. It is a ``centralized'' redistribution policy, which we call CR policy, in which part of the available stock at the outlets is sent to the warehouses and consolidated there before being sent back to their new destinations. No lateral transshipments are allowed.
	\item Figure~\ref{subfig:decentralized} represents a more standard (partially) decentralized approach, DR policy, that combines replenishment from the warehouses with lateral transshipments between outlets. There is no consolidation at the warehouses.
	\item Finally, Figure~\ref{subfig:general} represents the full generality of the policies captured by our model, GR policy, under which all previous possibilities are allowed.
\end{itemize}

Note that, although we assume in all models that decisions are made simultaneously for the whole network by a central planner, under CR policy the actual shipping can be completely monitored on-site from the warehouses. No such possibility exists when lateral transshipments are allowed. This centralized monitoring is particularly important for our partner, since some outlets are directly owned by the company and others are franchises.

As far as we know, CR policy has not been studied by past literature. This is not surprising since, intuitively, it seems inefficient as compared to DR and GR policies. Because of consolidation at the warehouses, order lead times get ``doubled'' under CR policy. Further, since some articles that might be transferred directly between outlets are sent first to the warehouse and then sent back to the destination outlet, transportation costs might ``double'' as well. Yet, as we discuss more deeply in Section~\ref{sec:policies}, CR policy also offers important advantages: redistribution becomes simpler and, since all shipments have a warehouse as the origin or the destination, the company has more bargaining power to negotiate the postage rates at their locations. This last observation is connected with an important feature of our analysis which, to the best of our knowledge, has not been dealt with by past literature: our optimal redistribution model accounts for the fact that articles are not sent individually, but in packages which may themselves be of different types. Thus, transportation costs depend on the number packages of each type sent between each pair of facilities. This implies that optimal redistributions will tend to fill packages as much as possible and, if possible, avoid lateral transshipments that involve a very small number of articles. This gives CR policy a potential advantage with respect to DR policy, since consolidation at the warehouse may result in a significant reduction on the total number of packages.\footnote{Think of a network with one warehouse and $n$ outlets. The total number of possible transshipments between outlets is $n^2$, whereas CR policy requires at most $2n$ movements: $n$ to send articles from outlets to the warehouse and $n$ to send them back the new destinations.} In Section~\ref{sec:policies} we present a numeric analysis in which we compare the optimal redistributions under CR, DR, and GR policies and discuss the main factors that drive the results. Further, the numeric analysis developed in Section~\ref{sec:TRP_num} shows that the model can be applied to obtain solutions for problems involving millions of integer variables.

The second main contribution of this paper is the design (Section~\ref{sec:TRP_def}) and analysis (Section~\ref{sec:TRP_num}) of a tailor-made algorithmic procedure to solve the resulting mathematical programming model. This is particularly relevant since, as usually in this field, we are confronted with a large-scale mixed-integer linear programming problem, which is hard to solve to optimality for instances arising in actual applications. The procedure we develop exploits the underlying network structure by first solving a relaxed version of the original problem and then recovering integrality by means of a suitably defined rounding scheme. We refer to this procedure as \emph{relaxed transferring-rounding-packing}, \TRPA, to represent each of its three phases. The performance of \TRPA is then benchmarked against the direct solution of the baseline models, \emph{transferring-packing}, which we call \TPA. Interestingly, although direct solution clearly dominates \TRPA for small instances, this domination gradually reduces as problem sizes increase, until eventually \TRPA outperforms \TPA.

\section{Contribution to Existing Literature}
Efficient inventory management is a hot research topic that has generated a lot of literature in the past and that continues to attract a lot of attention from the research community in the operations research's field. When attention is paid to stock redistribution in the presence of lateral transshipments, one important reference is the survey paper \cite{Paterson2011}, which includes a taxonomy of inventory models with lateral transshipments and take as the first layer in their classification the distinction between proactive and reactive redistribution models. Most research has focused on assessing the performance of different redistribution policies and, in particular, comparisons between proactive and reactive policies for transshipments are developed, for instance, in \cite{Hoadley1977}, \cite{Banerjee2003}, \cite{Burton2005}, and \cite{Lee2007}. Given its nature, proactive redistribution arises in the context of periodic review policies,\footnote{One exception is \cite{Paterson2012} where, in a continuous review setting, a transshipment prompted by a stockout is seen as an opportunity to do some proactive redistribution between the involved outlets.} and the issue of determining the timing of the redistribution is addressed in \cite{Agrawal2004}.

Additional complexity arises when redistributions may combine lateral transshipments between outlets with replenishments from warehouses. Our analysis belongs to this literature, already initiated in the seminal paper \cite{Gross1963} and later expanded in \cite{Karmarkar1977} for the single-period case and \cite{Karmarkar1981} and \cite{AboueeMehrizi2015} for the multi-period case. Contrary to these last papers, the assumption of negligible lead times is not present in \cite{Jonsson1987}, \cite{Diks1996,Diks1998}, and \cite{Tagaras2002}. In a recent paper, \cite{Firoozi2020} find that the combination of lateral transshipments and multiple sourcing in multi-echelon networks leads to substantial benefits.

While all of the above references consider different forms of proactive policies and possibly comparisons with reactive ones, a quick glance at Tables~2-4 in \cite{Paterson2011} shows that it is the literature on the purely reactive case the one that accounts for the majority of the research in the area. Most models consider the setting in which a period starts with the replenishment from the suppliers/warehouses, then demand is realized and, before demand is satisfied, it is possible to transship units between outlets (see, for instance, \cite{Hu2005}, \cite{Herer2006}, and \cite{Noham2014}). In the particular case of two-echelon supply chains in which decisions on transshipments and replenishments are taken simultaneously, \cite{Wee2005} study five different polices, depending on the allowed combinations for replenishments and lateral transshipments, and characterize the conditions under which each of them is optimal. Other papers studying the impact of lateral transshipments in connection with replenishment decisions are \cite{Bendoly2004}, \cite{Zhang2005}, and \cite{Tang2010}.

A significant part of the aforementioned literature studies multi-period settings in order to compare and understand different families of redistribution policies under uncertain demand. Each of the periods under study involves, typically in this order, replenishment decisions, demand realizations, and redistribution via lateral transshipments and/or back orders to fulfill demand. The complexity of these dynamic settings calls for important simplifications for the sake of tractability of the resulting models. This paper, on the other hand, instead of looking at the global inventory management process, develops a single-period model that focuses on the interim stage at which redistribution takes place. Once the object under study has been narrowed in this way, one can aim to build a model that can be readily applied to help in managerial decisions. In particular, the model we develop in Section~\ref{sec:model} allows for an arbitrary number warehouses, outlets, and items, referred to as stock keeping units (SKUs), and incorporates the following elements:

\begin{description}[style=nextline,topsep=0cm]
\item[Model's input:]
\begin{inlineitemize}
\anitem{Stock at warehouses.} \anitem{Stock at outlets.} \lastitem{Demand at outlets.}
\anitem{Types of packages.} \anitem{Shipping costs.} \lastitem{Shipping priorities.}
\anitem{Weight of each SKU.}
\end{inlineitemize}

Shipping costs may vary depending on the package type and on the facilities involved. Priorities are intended to give extra flexibility to the model, signaling which SKUs and/or outlets should be given more importance in the final redistribution.
\item[Model's output:]
\begin{inlineitemize}
	\lastitem{Units of each SKU sent between each pair of facilities.} 
	\lastitem{Packages of each type sent between each pair of facilities.}
\end{inlineitemize}
\end{description}

Past literature has associated a cost to each item to be sent or to each transshipment operation between two outlets, regardless of the number of items sent (see Tables 2-4 in \cite{Paterson2011}. We are not aware of any reference in which the packing component is incorporated into the cost structure of the redistribution process. This component is very important for our partner, since it performs its redistribution using package delivery companies and, hence, taking into account the specifics of their packages and pricing plans can have a substantial impact in the final solution.

\begin{figure}[!htbp]
\centering
\includegraphics[scale=0.4]{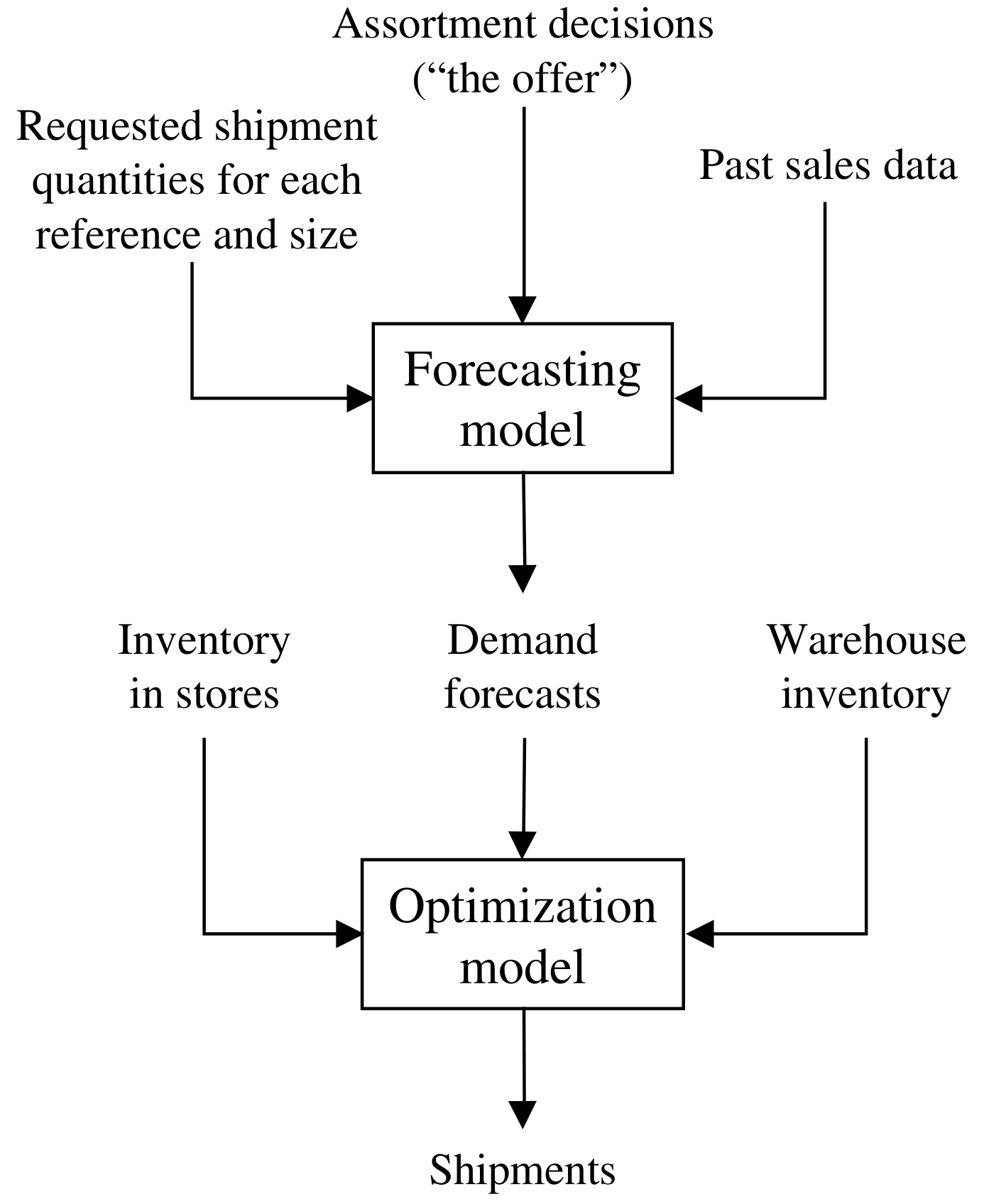}%
\caption{Scheme for Zara's periodic redistribution process as reported in \cite{Caro2010a} and \cite{Caro2010b}.}%
\label{fig:Zara}%
\end{figure}

Possibly the closest contribution to ours can be found in the analysis developed for Inditex Group, also in the clothing retail industry, in \cite{Caro2010a} and \cite{Caro2010b}.\footnote{We refer the reader to \cite{Wen2019} for a recent review on the use of Operations Research in fashion retail supply chain management.} Although they do not allow for lateral transshipments, they present a very detailed model for the interim redistribution optimization problem. As shown in Figure~\ref{fig:Zara}, taken from \cite{Caro2010a}, their optimization model uses essentially the same input as ours. Importantly, they also develop a detailed forecasting model to produce the estimates for future demand that are then fed into the optimization model.

\section{Model Discussion and Mathematical Formulation}
\label{sec:model}

We make the standard assumption that there is a central planner who is fully informed about the status of the supply chain and takes decisions with the objective of minimizing the system costs. More specifically, we take transportation costs as the main ingredient of the objective function and disregard holding costs. This approach is standard in the retail industry and, in particular, in the clothing industry since, as argued in \cite{Paterson2011}, handling costs are often dominant with respect to holding costs. This is one of the main features that makes proactive redistribution policies particularly suitable in this setting. Thus, although our model can be applied to compute optimal redistributions prompted by stockouts, it is primarily meant to be applied to (periodic) integrated network-wide redistributions.

\subsection{Main model: {\sc Transferring problem}}

The main objective is to develop an optimization model to redistribute stock across a network of facilities containing both warehouses and outlets. We refer to the associated optimization problem as the transferring problem, \tp, which we describe below. The main features of the model that distinguish it from previous literature are the following:

\begin{description}
	\item[F1.] The model can account for all forms of transfers of stock between facilities, including replenishments from the warehouses and lateral transshipments between outlets. Further, a novel contribution of this paper is to allow for the consolidation of outlet's excess stock at the warehouses before the replenishment from the warehouses is performed (see Figure~\ref{subfig:general}).
	\item[F2.] SKUs are sent in packages with limited capacity. Thus, transportation costs are not computed independently for each SKU, but depend on the packages sent between each pair of facilities. Therefore, the goal is to find a redistribution that optimizes the packages that travel through the network. 
	\item[F3.] Each outlet has associated two types of demand for each SKU: i) ``fixed demand'', which is known and must be met and ii) an estimate of ``variable demand'', which is not required to be met.
\end{description}

The next table describes all the elements involved in the definition of \tp:

%\begin{table}[!htbp]
\begin{center}
	\small
	%\begin{tabular}{m{5.5cm} m{9cm}}
	\begin{tabular}{ll}
	\hline
		\textbf{INPUT: Sets}&\\
		\hline \noalign{\vskip 0.1cm} 
		$\fac$ & Set of facilities (outlets and warehouses).\\
		$\ware\subset\fac$ & Set of warehouses.\\
		$\outl=\fac\setminus \ware$ & Set of outlets.\\
		$\sku$ & Set of stock keeping units.\\
		$\pack$ & Set of types of packages.\\
		$\mov \subseteq \fac \times \fac$ & Possible movements. {\footnotesize If $(i,j)\in \mov$, transfers from $i$ to $j$ are allowed.}\\
		& {\footnotesize Self-loops are not allowed: for each $i\in \fac$, $(i,i)\notin \mov$}.\\[0.2cm]
		\hline 
		\textbf{INPUT: Parameters}&\\
		\hline \noalign{\vskip 0.1cm} 
		$\sini_{is},\ i\in\fac,\ s\in\sku$ & Initial stock of SKU $s$ at facility $i$.\\
		$\fixd_{is},\ i\in\outl,\ s\in\sku$ & Fixed demand of SKU $s$ at outlet $i$ (has to be met).\\
		$\vard_{is},\ i\in\outl,\ s\in\sku$ & Variable demand of SKU $s$ at outlet $i$ (may not be met).\\
		$\pri_{is}\in[0,1],\ i \in \outl,\ s\in\sku$ & Priority given to outlet $i$ on SKU $s$.\\
		$\wei_s,\ s\in\sku$ & Weight of SKU $s$.\\
		$\capa_p,\ p\in\pack$ & Capacity of a package of type $p$.\\
		$\cost_{ijp},\ (i,j)\in\mov,\ p\in\pack$ & Cost of sending a package of type $p$ from facility $i$ to facility $j$. \\[0.2cm]
		\hline
		\textbf{OUTPUT: Variables}& \\
		\hline  \noalign{\vskip 0.1cm} 
		$\yvar_{ijp}\in \mathbb{Z}^{+},\ (i,j)\in\mov,\ p\in\pack$ & Number of packages of type $p$ sent from $i$ to $j$.\\
		$\xvar_{ijs}\in \mathbb{Z}^{+},\ (i,j)\in\mov,\ s\in\sku$ & Number of units of SKU $s$ sent from facility $i$ to facility $j$.\\[0.2cm]
		\hline
		\textbf{Aux. variables} {\scriptsize (to ease notation)}& \\
		\hline \noalign{\vskip 0.1cm} 
		%$\boldsymbol{Z}^{send}_{i,r}\in \mathbb{Z},\ i\in\fac,\ r\in\sku$ & $\sum_j X_{i,j,r}$. Sum of units of $r$ sent from store $i$. \prob{t} \\
		%$\boldsymbol{Z}^{rec}_{i,r}\in \mathbb{Z},\ i\in\fac,\ r\in\sku$ & $\sum_j X_{j,i,r}$. Sum of units of $r$ received by store $i$. \prob{t}\\
		$\sfin_{is}\in \mathbb{R},\ i\in\fac,\ s\in\sku$ & Final stock of SKU~$s$ in facility~$i$ defined as\\[1ex]
		\multicolumn{2}{c}{$\sfin_{is}=\sini_{is}+\sum_{j:(j,i)\in\mov}\xvar_{jis}-\sum_{j:(i,j)\in\mov}\xvar_{ijs}$.}  \\[0.2cm]
		\hline
		\textbf{Controllable coefficients} & \\
		\hline \noalign{\vskip 0.1cm} 
		$\alpha\in \mathbb R^+$ & Aggressiveness to meet variable demand. \\
		$\varepsilon \in \mathbb R^+$, $\varepsilon\approx 0$ & Technical parameter to break ties between solutions. \\
		\hline
	\end{tabular}
\end{center}
	%\caption{Sets, parameters, and variables involved in the transferring problem.}
	%\label{tab:spv}
%\end{table}

We are now ready to present and discuss the formulation of \tp:
\begin{align}
	\multispan2{\sc{Transferring Problem: \tp} \hfil} \notag \\
	% Obj1: Shipping Cost
	\min_{\xvar,\yvar,\sfin} & \sum_{\substack{(i,j)\in\mov \\ p\in\pack}} \cost_{ijp} \yvar_{ijp} \quad \label{cons:T_obj1} \\
	% Obj2: Penalty for not satisfying variable demand
	&\qquad + \alpha \sum_{\substack{i\in\fac\\s \in \sku}} \pri_{is} \max\{0, \fixd_{is} + \vard_{is} - \sfin_{is}\}  \label{cons:T_obj2} \\
	% Obj3: Penalty for unnecessary movements
	&\qquad + \varepsilon \sum_{\substack{(i,j)\in\mov\\s \in \sku}} \xvar_{ijs} \label{cons:T_obj3} \\
	% Cons1: Formula for final stock
	\text{s.t.} \qquad & \sfin_{is}=\sini_{is}+\sum_{j:(j,i)\in\mov}\xvar_{jis}-\sum_{j:(i,j)\in\mov}\xvar_{ijs}, \quad i \in \fac,\, s \in \sku \label{cons:sfin}\\
	%
	% Cons2: Fix Demand has to be satisfied
	& \sfin_{is}  \geq \fixd_{is}, \quad i\in\outl,\,s \in \sku\label{cons:T_minitem}\\
	%
	% Cons 3: At most final stock \fixdand + \vardand
	%& \sfin_{is} \leq \fixd_{is}+\vard_{is}, \quad i\in\outl,\ s \in \sku  \label{cons:T_maxitem}\\
	%
	% Cons 4: Don't send if you need
	& \sum_{\substack{j:}{(i,j)\in\mov}}\xvar_{ijs} \leq  \max\{0,\, \sini_{is} - \fixd_{is}\}, \quad i\in\outl,\, s \in \sku \label{cons:T_dsiyn}\\
	%
	% Cons 5: Capacity constraints
	& \sum_{s \in \sku} \wei_s \xvar_{ijs}\leq \sum_{p\in\pack} \capa_p\yvar_{ijp}, \quad (i,j)\in\mov \label{cons:T_leqcap} \\
	& \xvar_{ijs} \in {\mathbb{Z}^{+}},\quad (i,j)\in\mov,\, s \in \sku  \notag \\
	& \yvar_{ijp} \in \mathbb{Z^{+}}, \quad (i,j)\in\mov,\, p\in\pack. \notag 
\end{align}

The \textbf{objective function} has two main terms. First, \eqref{cons:T_obj1} represents the overall transportation costs associated with the stock redistribution proposed by the model. Second, in \eqref{cons:T_obj2} we have a term that penalizes unmet variable demand. Coefficient $\alpha$ allows to control how aggressive to be in satisfying variable demand, \emph{i.e.}, to control the level of proactivity (and riskiness) sought by the model. In particular, one can obtain purely reactive redistributions by setting $\alpha=0$. Finally, there is a third ``technical'' term, \eqref{cons:T_obj3}, meant to break ties and reduce the set of optimal solutions (which should improve solver's performance). More specifically, given two solutions that result in the same packages being sent and the same final stocks, one should give priority to the redistribution that involves less transfers.\footnote{One could have, for instance, an optimal solution in which the same SKU is sent from $i$ to $j$ and from $j$ to~$i$, just because there is available space in the packages that must be sent between the two facilities. Thus, this ``useless'' movement entails no cost and could be optimal if $\varepsilon=0$ in \eqref{cons:T_obj3}.} Given the role of this third term, the value of $\varepsilon$ has to be close to zero.

We turn now to discuss the \textbf{constraints}: \eqref{cons:sfin} just defines the auxiliary variables for the final stock, $\sfin$, to ease notation in the rest of the model, \eqref{cons:T_minitem} ensures that fixed demand is met, \eqref{cons:T_dsiyn} states that an outlet can only send items that it has in excess after meeting its own fixed demand. Finally, \eqref{cons:T_leqcap} represents the linking constraints for the $\xvar$ and $\yvar$ variables, requiring that the total number of items shipped between two facilities cannot exceed the total capacity of the packages sent between them.

Note that constraints in~\eqref{cons:T_dsiyn} do not apply to warehouses which, in particular, are even allowed to send stock ``before'' actually having it available. This is only relevant when set $\mov$ allows for transfers from the outlets to the warehouses. This stock can be consolidated at the warehouses and then sent back to other outlets. This feature of the model is essential to capture as a particular case different variants of the so called centralized redistribution schemes, in which lateral transshipments are not possible. 

\begin{figure}[!htbp]
\centering
\begin{tikzpicture}[scale=1,>={Latex[length=2mm,width=3mm]},font=\footnotesize,every node/.style={font=\footnotesize},estilo/.style={rectangle,draw, align=center, inner sep=0.2cm,fill=gray!50, rounded corners}]

\node[estilo] (w) at (0,0) {Warehouse};
\node[estilo] (o1) at (3,1.5) {Outlet 1 (O1)};
\node[estilo] (o2) at (3,-1.5) {Outlet 2 (O2)};

\draw[->>,very thick] (w) -- (o1);
\draw[->>,very thick] (o1) -- (o2);
\draw[->>,very thick, dashed, gray] (w) -- (o2);

\node[anchor=west] (table) at (5,0){
	\begin{tabular}{cccc}
	 & Warehouse & Outlet 1 & Outlet 2 \\
		\hline
		$\sini_{is_1}$ & 5 & 0 & 0\\
		$\sini_{is_2}$ & 0 & 1 & 0\\
		$\sini_{is_3}$ & 5 & 1 & 0\\
		\hline
		$\fixd_{is_1}$ & 0 & 1 & 0\\
		$\fixd_{is_2}$ & 0 & 0 & 1\\
		$\fixd_{is_3}$ & 0 & 1 & 1\\
		\hline
	\end{tabular}
};
\end{tikzpicture}
\caption{Illustrating the role of constraints~\eqref{cons:T_dsiyn}.}
\label{fig:dsiyn}
\end{figure}
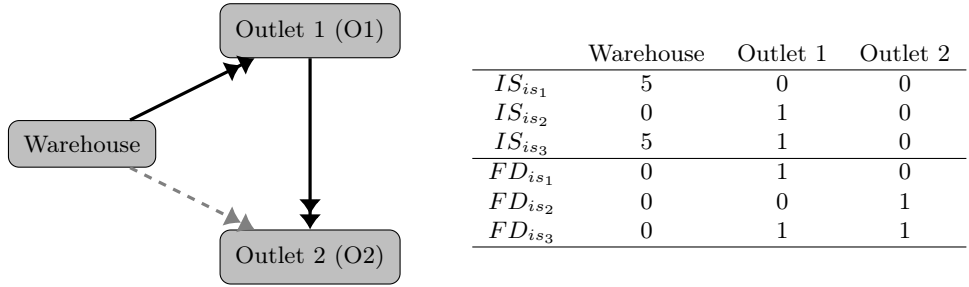

Finally, it is worth noting that, depending on the specific setting under study, it might be convenient to weaken constraints in~\eqref{cons:T_dsiyn} and just require that
\begin{equation}
\sum_{\substack{j:}{(i,j)\in\mov}}\xvar_{ijs} \leq  \sini_{is}, \quad i\in\outl,\, s \in \sku. \label{cons:T_sendifav}
\end{equation}
This would allow redistributions in which an outlet sends units of a SKU that it needs to meet its fixed demand because it will receive them back from some other facilities. The small network in Figure~\ref{fig:dsiyn} illustrates why this possibility may reduce the total number of packages used for redistribution.

Consider the initial stocks and fixed demands in the table of Figure~\ref{fig:dsiyn}, and assume that shipments from the outlets to the warehouse are not allowed. Then, in any feasible solution of the associated redistribution problem the warehouse has to send a package to O1 with one unit of $s1$ and O1 has send a package to O2 with one unit of $s2$. Further, under constraints~\eqref{cons:T_dsiyn}, an additional package must be sent from the warehouse to O2 with one unit of $s_3$. However, if~\eqref{cons:T_dsiyn} is weakened to~\eqref{cons:T_sendifav}, this additional package is not needed, since one can send one unit of $s3$ from O1 to O2 and another one from the warehouse to O1 using the already available packages (provided they have enough capacity).\footnote{The complexity of solving \tp does not change depending on whether~\eqref{cons:T_dsiyn} or~\eqref{cons:T_sendifav} is used. Further, although hereafter we restrict attention to constraints in~\eqref{cons:T_dsiyn}, all the analysis presented in this paper has been replicated replacing them with those in~\eqref{cons:T_sendifav} and the results are qualitatively very similar.}

\subsection{Model features: Discussion}
\label{sec:features}

As we have already mentioned, the main features of this model that distinguish model \tp from previous literature are \textbf{F1}, \textbf{F2}, and \textbf{F3}:

\begin{itemize}
	\item \textbf{F1} allows to bring different redistribution policies under the same umbrella (see Figure~\ref{fig:redist_schemes}) and eases the comparisons between them, as we illustrate in Section~\ref{sec:policies}. Different specifications of set $\mov$ in the formulation of \tp correspond with different policies. As we have already mentioned, a special case is the policy used by our partner: a centralized redistribution scheme without lateral transshipments in which outlets send part of their stocks to the warehouses, where they get consolidated with the stocks available there and used for the replenishment. Note that this consolidation implies that the lead times get doubled, since warehouses must wait for the items sent by the outlets before proceeding with the replenishment. Although this increase in lead times is not captured explicitly by our static model, it can be indirectly accounted for by appropriately increasing the costs of the packages sent to warehouses.
	
	\item \textbf{F2} adds the packages' layer to the model, which brings it closer to real applications and which, as far as we know, was absent in previous literature. In return, we get a much more complex optimization model. Indeed, note that if we take away the packing component from the model by removing constraints in~\eqref{cons:T_leqcap} and replacing $\yvar_{ijp}$ with $\xvar_{ijs}$ in the objective function, we get a network flow problem that can be decomposed in independent flow problems, one for each $s\in \sku$. Not only each of these problems would be much smaller than the original problem but, because of the unimodularity property underlying these flow problems, they could essentially be solved as linear problems.\footnote{Because of this underlying structure, when solving	\tp one might greatly benefit from specialized solution methods such as Benders' decomposition or branch and price algorithms, as we briefly discuss in Section~\ref{sec:TPA}.}
	\item \textbf{F3} concerns the explicit specification of fixed and variable demand, which allows to account for both reactive and proactive redistribution policies and for intermediate policies as well.
\end{itemize}

The role of coefficient $\alpha$ in the objective function is similar to the role of parameter $K$ in \cite{Caro2010a} and \cite{Caro2010b}, where the authors argue that ``\emph{value $K$ can thus be interpreted as the unit opportunity cost of shipping inventory to a particular store and is meant as a control lever allowing the model user to affect its output}''. They then go on to argue that this parameter allows to control how ``conservative'' or ``aggressive'' is the solution proposed by the model. In particular, one may choose to be more or less aggressive depending on the reliability of the estimates for the variable demand. Interestingly, in our model one can use the $\pri_{is}$ parameters to make this aggressiveness specific to pair facility-SKU, $(i,s)$.

Ideally, \tp model should be part of a more complex framework, including a first stage for the forecasting of future demand. In particular, the estimates for the variable demand might i)~be made the different outlets in a decentralized way,\footnote{This is how our partner operated at the time in which this project was developed.} ii)~come from a centralized forecasting model like the one developed in \cite{Caro2010a} and \cite{Caro2010b}, or iii)~be a combination of both types of estimates.

Problem \tp is a linear optimization problem with integer variables, which makes it an NP-hard problem. Since realistic applications may involve a large number of facilities and a very large number of SKUs, solving \tp to optimality is a very challenging task. Sections~\ref{sec:TRP_def} and~\ref{sec:TRP_num} are devoted to present a tailor made approximated procedure to solve \tp and to compare its performance with the direct solution of the underlying model.

\subsection{Enhancing the final solution: \sc{Packing problem}}\label{subsec:packingdef}

The packing component of the transferring problem defined above has one limitation that is worth discussing. Constraints in~\eqref{cons:T_leqcap} require that the total weight of the items send from $i$ to $j$ does not exceed the capacity of the packages sent from $i$ to $j$. Yet, this would allow, for instance, to send three items of weight 3 in two packages of capacity 5. However, this is not feasible without ``dividing'' one of the three items, which is not possible and, therefore, a third package would be needed.

Addressing the above problem requires to explicitly model a \emph{packing problem} between each pair of facilities.\footnote{Refer, for instance, to Chapter~18 in \cite{Bernhard2006} and the survey papers \cite{Lodi2002}, \cite{Coffman2013}, and \cite{Christensen2017}.} These problems are themselves difficult mathematical programming problems (NP-hard) and incorporating them into problem \tp would result in an unmanageable optimization problem. Despite this limitation it is reasonable to expect that, in most practical applications, the impact of situations like the one described above on the overall transportation costs will be relatively small.\footnote{Since the packing problem does not affect the $\xvar_{ijs}$ variables, its impact on the objective function of problem \tp is only in term~\eqref{cons:T_obj1}. Table~\ref{tab:packages} in Section~\ref{sec:results_instances_sml} shows that, in the different test sets under study, the average increase in this term after \pp is run varies from $10$ to $18\%$ (the impact on the whole objective function is smaller and varies substantially depending on the value of $\alpha$).} Thus, instead of incorporating the packing problems into problem \tp, we solve them at a post-processing phase on \tp's final outcome in order to provide the final user with detailed information about the composition of each package to be sent.

More precisely, once \tp is solved, we have a redistribution that specifies, for each pair of facilities, the number of items of each SKU ($\xvar_{ijs}$) and the number of packages ($\yvar_{ijp}$) to be sent between them. In order to determine the specific composition of each package and to see whether or not some additional package might be needed, we solve a packing problem, \pp, for each pair of facilities. The next table describes all the elements involved in the definition of problem \ppm{ij}, with $(i,j) \in \mov$:
%\begin{table}[!h]
\begin{center}
	\small
	%\begin{tabular}{m{5.5cm} m{9cm}}
		\begin{tabular}{ll}
	\hline
		\textbf{INPUT: Sets}&\\
		\hline \noalign{\vskip 0.1cm} 
		$\pack$ & Set of types of packages.\\
		$\sku$ & Set of stock keeping units.\\[0.2cm]
		\hline 
		\textbf{INPUT: Parameters}&\\
		\hline \noalign{\vskip 0.1cm} 
		$\avap_p$, $p\in \pack$ & Number of available packages of type $p$.\\
		$\capa_p$, $p\in\pack$ & Capacity of a package of type $p$.\\
		$\cost_{p}:=\cost_{ijp}$, $p\in\pack$ & Cost of sending a package of type $p$. \\
		$T_s:=\xvar_{ijs}$, $s \in \sku$ & Number of units to transfer of SKU $s$. \\
		$\wei_s$, $s\in\sku$ & Weight of SKU $s$.\\[0.2cm]
		\hline \noalign{\vskip 0.1cm} 
		\textbf{OUTPUT: Variables}& \\
		\hline \noalign{\vskip 0.1cm} 
		$\xpack_{pks}\in \mathbb{Z^{+}}$, $p\in \pack, k\in \{1,\ldots, \avap_p\}, s\in \sku$  & Number of units of SKU $s$ sent in package $k$ of type $p$. \\
		$\ypack_{pk}\in \{0,1\}$, $p\in \pack, k\in \{1,\ldots, \avap_p\}$ & Represents whether or not package $k$ of type $p$ is sent.\\[0.2cm]
		\hline
	\end{tabular}
\end{center}

We now present the formulation of the packing problems \ppm{ij}, one for each $(i,j) \in \mov$:
\begin{align}
	\multispan2{{\sc Packing Problem:} \ppm{ij}, with $(i,j) \in \mov$ \hfil} \notag \\[0.1cm]
	\min_{\xpack,\ypack} & \sum_{p\in \pack} \sum_{k=1}^{\avap_p} \cost_{p}\boldsymbol{y}_{pk} \notag \\[0.1cm]
	\text{s.t.} \qquad & \sum_{s \in \sku} \wei_s \xpack_{pks}\leq \capa_{p} \ypack_{pk}, \quad p \in \pack,\,k \in \{1,\ldots,\avap_p\} \notag \\
& \sum_{p\in \pack} \sum_{k=1}^{\avap_p} \xpack_{pks}=T_s, \quad s\in \sku  \notag \\
& \xpack_{pks} \in \mathbb{Z^{+}}, \quad p\in \pack,\, k\in \{1,\ldots, \avap_p\},\, s\in \sku \notag \\ 
& \ypack_{pk} \in \{0,1\},\quad p\in \pack,\, k\in \{1,\ldots, \avap_p\}.  \notag 
\end{align}	

Once all these problems have been solved we would have the exact number of packages needed to implement the stock redistribution associated with the $\xvar_{ijs}$ variables obtained when solving \tp. 

\section{Comparison of Redistribution Policies}
\label{sec:policies}

\subsection{Redistribution policies}

One of the main drivers of this research was that our partner wanted to compare the centralized redistribution policy in use by the company with alternative policies. In this section we illustrate how the versatility of the \tp formulation can be used to develop such a comparison.\footnote{The results reported in this section correspond with the direct solution of \tp problem, but qualitatively identical results were obtained after enhancing the solutions by running \pp on the final outcomes.} More precisely, we focus on the following three redistribution policies, already outlined in the introduction:
\begin{description}
    \item[Centralized redistribution, CR.] All the redistribution process is centralized at the warehouses (Figure~\ref{subfig:centralized}). Lateral transshipments are not allowed, but outlets can still send stock to each other via the warehouses. This was the policy in place at our partner at the time of this study. This is achieved by defining $\mov$ to be the set obtained after removing from $\fac\times\fac$ all the elements in $\outl\times\outl$ and the self-loops.
    \item[Decentralized redistribution, DR.] Combines replenishment from the warehouses with lateral transshipments between outlets (Figure~\ref{subfig:decentralized}). There is no consolidation at the warehouses. Set $\mov$ is now defined by removing from $\fac\times\fac$ all the elements in $\outl \times \ware$ and the self-loops.
    \item[General redistribution, GR.] This is the most general policy where $\mov$ is obtained by removing the self-loops from $\fac\times \fac$ (Figure~\ref{subfig:general}). 
\end{description}

In the case of our partner, Nanos, there were three main reasons behind the use of CR policy: i) given the large volume of shipments involving the warehouses, the company has significantly more bargaining power to negotiate with the different package delivery companies the postage rates of these shipments, ii) all redistribution can be completely monitored from the warehouses, which is particularly important given that some of Nanos' outlets are franchises not owned by the company, and iii) given that Nanos is a high-end clothing retail company, the number of items usually involved in these redistributions is not as large as for other middle and low-end retail companies. Hence, in situations where DR might call for shipments in which some packages are sent almost empty just to satisfy a very specific demand of a given outlet, consolidation at the warehouses may end up reducing the total number of packages.

In Section~\ref{sec:numeric_policies} we briefly explain the methodology used to generate the test set on which the above three policies are compared and present the numeric results. Note that, since GR includes CR and DR as special cases, GR will always deliver the best solutions according to model \tp. The idea is to build upon point i) above and see how CR and DR separate from GR policy as the costs of the transfers involving the warehouses vary relative to the costs of transfers between outlets (lateral transshipments).

%\subsection{Generation of the test set}
%\label{sec:instace_generation}

\subsection{Numeric comparison}
\label{sec:numeric_policies}

All the test sets used in this paper, both in this section and in Section~\ref{sec:TRP_num}, have been randomly generated using an instance simulator implemented in Python, and whose rules ensure feasibility of the resulting instances.\footnote{A detailed explanation of these rules, including the pseudocodes of the basic procedures, can be found in Appendix~\ref{sec:appendix_generator}.} This simulator allows to control different parameters during the generation process, which include the number of elements of $\ware$, $\outl$, $\pack$, $\sku$, and the total stock of the network (sum of the number of units of all SKUs).

For the comparative analysis of the three redistribution policies, CR, DR, and GR, we have generated 50 instances using the aforementioned simulator. All these instances were of the same size, with $\card{\ware}=1$, $\card{\outl}=10$, $\card{\pack}=2$, $\card{\sku}=10$, and a total of $1000$ items involved in the redistribution process. This results in instances with 1430 variables (1320 integer) and where the number of constraints depends on the specific redistribution policy. In particular, the largest number corresponds with GR policy, with 430 constraints. All 50 instances of problem \tp were solved to optimality in around 5 minutes.

Different variations of each of the 50 instances generated above are solved in order to develop the comparison. First, each instance was solved for four different values of the aggressiveness factor, $\alpha \in \{0, 0.1, 10, 100\}$, resulting already in 200 instances. Then, to represent different balances of the shipping rates of packages involving the warehouses with respect to direct shipments between outlets, a scaling factor was introduced to multiply the former. More precisely, each of the above 200 instances was solved for 13 different values of this scaling factor: 0.01, 0.1, 0.25, 0.5, 0.75, 0.9, 1, 1.1, 1.5, 2, 5, 10 and 100. Values smaller than 1 lead to instances in which there is a ``discount'' in shipments involving the warehouses and values greater than 1 correspond to ``penalizations'' in those shipments. Since each instance is solved for the three redistribution policies by varying the definition of the set $\mov$, the total number of instances solved in this comparative analysis adds up to $200\times13\times 3 =7800$, 2600 for each redistribution policy.

\begin{figure}[!htbp]
\centering
\begin{subfigure}{.48\textwidth}
  \centering 
	\includegraphics[height=5.35cm,viewport=14mm 8mm 242mm 165mm,scale=0.5,clip]{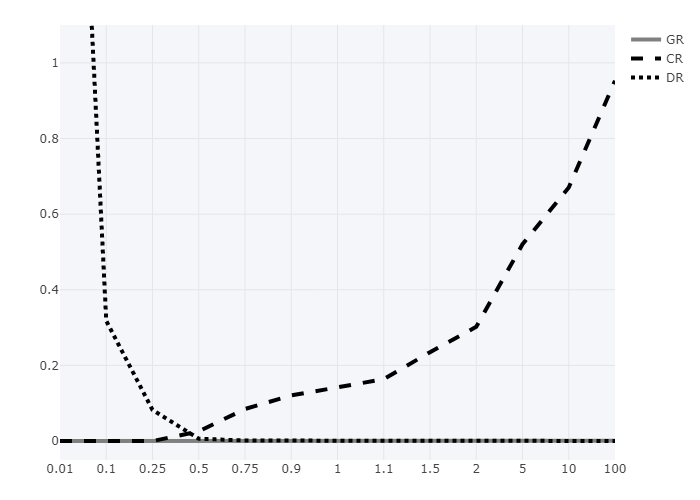}
	\caption{Objective function of \tp.}
	\label{fig:comp1}
\end{subfigure}
\hspace{-0.2cm}
\begin{subfigure}{.48\textwidth}
	\centering
	\includegraphics[height=5.35cm,viewport=14mm 8mm 222mm 165mm,scale=0.5,clip]{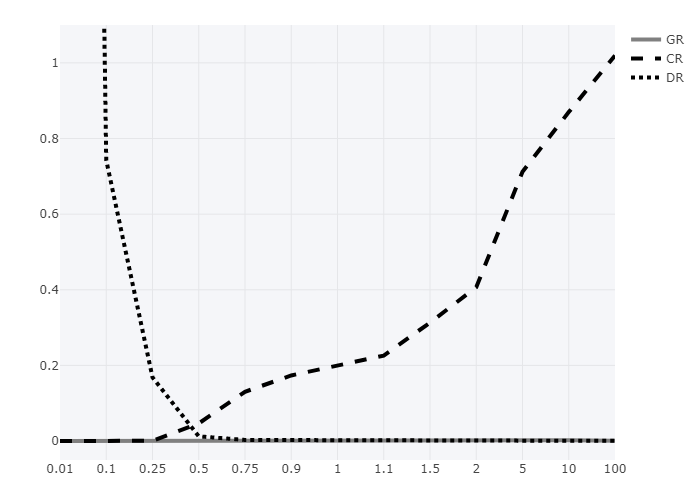}
	\caption{Transportation costs in \tp.}
	\label{fig:comp2}
\end{subfigure}
\caption{Average worsening with respect to best redistribution policy as a function of the scaling factor for the shipping cost of packages involving warehouses.}
\label{fig:shiping_cost_centralized_vs_decentralized}
\end{figure}

Figure~\ref{fig:shiping_cost_centralized_vs_decentralized} summarizes the results. For each value of the scaling factor, the numbers represented in Figure~\ref{fig:comp1} are computed as follows:
\begin{itemize}
	\item For each instance let $objCR$, $objDR$, and $objGR$ be the objective function at optimality of the different policies. Let $best$ denote the minimum of the three values. Then, the relative performance of each policy with respect to the best one is computed as follows:
	\[
	\frac{objCR-best}{best}, \quad \frac{objDR-best}{best}, \quad \frac{objGR-best}{best}.
	\]
	These numbers represent how many times more costly each policy is with respect to the best one.
	\item For each value of the scaling factor and each redistribution policy, we compute the average performance in the corresponding 200 instances. 
	\item We now have a unique number for each value of the scaling factor and each redistribution policy, which we represent in Figure~\ref{fig:comp1} (with the value of the scaling factor in the $x$ axis).
\end{itemize} 
Figure~\ref{fig:comp2} contains a similar representation in which only the transportation costs at optimality are considered, (term \eqref{cons:T_obj1} in the objective function).

As expected, GR policy is always optimal, since anything feasible for either CR or DR is also feasible for GR. Then, Figure~\ref{fig:shiping_cost_centralized_vs_decentralized} illustrates the behavior of CR and DR. For small values of the scaling factor, which correspond with reduced rates for shipments involving the warehouses, CR policy is optimal, \emph{i.e.}, lateral transshipments are too expensive relative to the consolidation at the warehouses. Conversely, for large values of the scaling factor, it is DR policy the one that becomes optimal. Interestingly, when the scaling factor takes value 1, \emph{i.e.}, all postage rates are generated using the same distribution and no scaling is performed, DR policy is again as good as GR policy. Looking at the figure, it seems that the turning point is around $0.5$, meaning that a $50\%$ discount on the postage rates involving the warehouses is needed to make CR policy competitive with DR policy. Note that it is precisely for those values where the two policies exhibit a similar performance where GR policy strictly improves upon both of them. This is because GR benefits from the combined use of the movements allowed by CR and DR.

The discussion above should just be seen as a conceptual qualitative exercise to illustrate the versatility of \tp model, and quantitative conclusions should note be drawn given that the analysis is based on randomly generated data. However, we believe that Figure~\ref{fig:shiping_cost_centralized_vs_decentralized} highlights the importance of designing suitable decision support systems. They should be aimed at helping company managers to choose between different redistribution policies depending on the specifics of the network-wide inventory when such a redistribution is called upon.

%A la vista de los resultados parece que una política centralizada únicamente debería de tener sentido si el precio del coste de los paquetes que pasan por almacén tiene un descuento mayor del 75\% con respecto al precio de un envío ordinario. Por contra, una política descentralizada tendría sentido únicamente en el caso de que el costes de los paquetes que pasan por almacén tuviesen un descuento menor del 50\%. Este análisis pone de manifiesto que la política centralizada no tiene sentido si el coste de los paquetes que pasan por almacén no tiene un descuento superior al 50\%, ya que para realizar un envío entre dos tiendas será necesario realizar dos envíos (tienda-almacén, almacén-tienda) en vez de uno (tienda-tienda). Por último, la tercera política es la clara ganadora ya que permite adaptarse al los precios de los envíos para tomar la mejor decisión. 

%Estos resultados están condicionados a la generación de las instancias, pero en caso de tener que decantarse por una política de redistribución, queda patente la necesidad de estudiar las condiciones en los precios de los paquetes para no incurrir en costes muy elevado al decidirse por una política no beneficiosa.

\section{\TPA: The Direct Solution Procedure}\label{sec:TPA}
Recall that, as we discussed in Section~\ref{subsec:packingdef}, incorporating into problem \tp the constraints associated to the \ppm{ij} packing problems, with $(i,j) \in \mov$, would result in an unmanageable optimization problem. Instead, since one would not expect these additional packing constraints to have a large impact on the final cost, we proposed the two-phase procedure \tppp, which we call \TPA. This procedure delivers an approximate solution of the overall redistribution problem: the cost of the optimal solution of \tp provides a lower bound on the overall optimal redistribution cost and the solution of \tppp provides an upper bound. Thus, the gap between these two bounds provides information about how large the approximation error of \TPA might be for each given instance. 

The core and computationally most demanding phase of \TPA is the solution of problem \tp. The redistribution problems meant to be solved by \tp quickly grow in size with the number of facilities and SKUs. Thus, the computational requirements for the exact solution of \tp of real size instances may become prohibitive and the design of algorithms for its efficient solution is also a matter of high relevance. Recall that in Section~\ref{sec:features}, when discussing the features of \tp, we mentioned that it has an underlying network structure (essentially a transportation problem), which is ``destroyed'' by the packing component. Without the latter, the problem could be decomposed in independent network flow problems, one for each $s\in \sku$. Hence, a natural approach to solve \tp would be to implement some decomposition technique that can exploit the underlying network structure such as: i) a branch and price algorithm to handle the complicating constraints on the packing variables and get a decomposition structure with independent minimum cost flow subproblems for each SKU,\footnote{Refer, for instance, to \cite{Vanderbeck1996}, \cite{Barnhart1998}, and \cite{Vanderbeck2000,Vanderbeck2011}.} ii) a Benders algorithm in which the complicating packing variables are handled by the master problem, obtaining again independent minimum cost flow subproblems for each SKU,\footnote{Refer, for instance, to \cite{Benders1962}, \cite{Geoffrion1972}, and \cite{Rahmaniani2017}.} or iii)~some specialization of Lagrangian relaxation, penalizing the complicating constraints in the objective function.\footnote{Refer, for instance, to \cite{Conejo2006}.} 

Instead of going over the exercise of studying the efficiency of the above classical approaches in our current problem, we have developed a new algorithmic idea, which we discuss in Section~\ref{sec:TRP_def} below. 

\section{\TRPA: An Algorithm for Large-Scale Problems}
\label{sec:TRP_def}

The direct solution of \tppp behind \TPA is replaced with the solution of a relaxed version of \tp, \rtp, combined with a tailor made rounding scheme, \rp, to recover a feasible solution of \tp that is then passed to \pp. We refer to this \rtprppp procedure as \TRPA. In this section we discuss the three phases of the \TRPA algorithm and in Section~\ref{sec:TRP_num} we present the numeric results to illustrate the following features: i) \TRPA delivers good solutions, comparable to those obtained by \TPA and, at the same time, ii) it provides better scaling for large instances. Importantly, \TRPA  should not only be seen as an alternative to \TPA for large instances, but its underlying principles can be used in conjunction with the decomposition techniques mentioned in Section~\ref{sec:TPA}, since the rounding scheme in \rp might also be used to speed up the (approximate) solution of the master and/or subproblems associated to these techniques.

\subsection{\TRPA. Phase 1: Relaxed transferring}\label{subsec:relaxed}

Relaxing the integrality of the $\yvar$ variables in \tp may have a big impact both on the final cost and on the underlying redistribution. On the other hand, a relaxation of the $\xvar$ variables might just have a second order impact in the cost function and lead to redistributions relatively close to those that one would get without relaxations. Further, such a relaxation might notably reduce solution times with respect to problem \tp, specially for large instances. With this idea in mind we define problem \rtp, which has just two differences with respect to problem \tp:
\begin{enumerate}
    \item The $\xvar_{ijs}$ variables can take continuous values. As will be clear from the rounding scheme discussed in Section~\ref{subsec:rounding}, one of the main limitations of the \rtprp procedure is that additional packages may be needed after rounding the continuous $\xvar_{ijs}$ variables produced by problem \rtp. The second modification of problem \tp aims at mitigating this limitation.
    \item There is a new parameter, $\delta \in (0,1]$, which is added to the capacity constraints~\eqref{cons:T_leqcap}, so that they become
    \[\sum_{s \in \sku} \wei_s \xvar_{ijs}\leq\sum_{p\in\pack} \delta\, \capa_p\yvar_{ijp}, \quad (i,j)\in\mov. \]
		This parameter controls how much the packages can be filled in problem \rtp. A value $\delta=0.8$ means that only $80\%$ of the available capacity can be used. This reserved capacity may be useful in the rounding stage, by reducing the number of packages that must be added in it.
\end{enumerate}

When choosing the value of $\delta$ there is trade-off that must be taken into account, since, the smaller $\delta$ is, the larger the number of packages needed in problem \rtp.\footnote{Our numerical analysis suggests that $\delta$ should not be smaller than $0.8$, although the precise value depends on the specifics of each application.}

\subsection{\TRPA. Phase 2: Rounding}\label{subsec:rounding}

Once a solution of problem \rtp has been obtained, the next step is to find rounded values of the $\xvar$ variables that is feasible to problem \tp. In order to accomplish this, we develop a rounding scheme, \rp, that builds upon the classic matrix rounding problem, introduced in \cite{Bacharach1966}.\footnote{See also \cite{Cox1982}, \cite{Doerr2006}, and \cite{Salazar2006}.} Computationally, the rounding scheme \rp is much less demanding than both \tp and \rtp because of the following two features: i) it can be seen as a special case of the minimum cost network flow problem,\footnote{Refer, for instance, to \cite{Ahuja1993}.} so its underlying unimodularity allows to solve it directly as a linear problem, ignoring the integrality constraints and ii) it is separable in independent problems for the different SKUs, which speeds up solution times even further. The next table describes all the elements involved in the definition of problem \rpm{s}, with $s\in \sku$:
%\begin{table}[!h]
\begin{center}
	\small
	%\begin{tabular}{m{5.5cm} m{9cm}}
		\begin{tabular}{ll}
	\hline
		\textbf{INPUT: Sets}&\\
		\hline \noalign{\vskip 0.1cm} 
		$\mov \subseteq \fac \times \fac$ & Possible movements.\\[0.2cm]
		\hline 
		\textbf{INPUT: Parameters}&\\
		\hline \noalign{\vskip 0.1cm} 
		$X^{Rel}_{ijs} \in \mathbb{R^{+}}$, $(i,j)\in\mov$ & Number of units of SKU $s$ sent from facility $i$ to facility $j$. \\
		& {\scriptsize (in the solution whose rounding is needed).}\\
$Z^S_{is}:=\displaystyle\sum_{\substack{j:(i,j) \in \mov}} X^{Rel}_{ijs} \in \mathbb{R^{+}}, (i,j)\in \mov$  & Number of units of SKU $s$ sent from facility $i$. \\[0.5cm]
$Z^R_{is}:=\displaystyle \sum_{\substack{j:(j,i) \in \mov}} X^{Rel}_{jis} \in \mathbb{R^{+}}, (i,j)\in \mov$ & Number of units of SKU $s$ received by facility $i$.\\[0.5cm]
$B_{is}:=Z^R_{is}-Z^S_{is} \in \mathbb{R}, i \in \fac$ & Balance between units send and received at facility $i$. \\
$\hat{c}_{ijs}\in \mathbb{R}$, $(i,j)\in\mov$ & Cost associated to each unit of SKU $s$ sent from $i$ to $j$. \\[0.2cm]
		\hline \noalign{\vskip 0.1cm} 
		\textbf{OUTPUT: Variables}& \\
		\hline \noalign{\vskip 0.1cm} 
		$\xvar{ijs}\in \mathbb{Z^{+}}$, $(i,j) \in \mov$  &  Number of units of SKU $s$ sent from facility $i$ to facility $j$. \\[0.2cm]
		\hline
	\end{tabular}
\end{center}

We are now ready to present and discuss the formulation of problem \rpm{s}, with $s\in \sku$:
\begin{align}
& \text{{\sc Rounding Problem:} \rpm{s}, with $s\in \sku$} \notag \\[0.1cm]
\min_{\xvar} & \sum_{\substack{(i,j)\in\mov}}  \hat{c}_{ijs} \xvar_{ijs}\\
%
% Box constraints for each variable
\text{s.t. }\quad & \Big\lfloor X^{Rel}_{ijs}\Big\rfloor \leq \xvar_{ijs} \leq \Big\lceil X^{Rel}_{ijs}\Big\rceil, \quad (i,j) \in \mov \label{eq:round_fac}\\[0cm]
%
% UB-LB Columns
&\Big\lfloor Z^S_{is} \Big\rfloor \leq \sum_{\substack{j:(i,j) \in \mov}}\xvar_{ijs} \leq \Big\lceil Z^S_{is}\Big\rceil, \quad i\in\fac \label{eq:round_bal_fac}\\[0cm]
%
% UB-LB Rows
&\Big\lfloor Z^R_{is}\Big\rfloor \leq \sum_{\substack{j:(j,i) \in \mov}}\xvar_{jis} \leq \Big\lceil Z^R_{is}\Big\rceil,\quad i \in\fac \label{eq:round_fac_bal}\\[0cm]
%
% UB-LB Prods
%&\bigg\lfloor \sum_{s \in \sku} X^{Rel}_{ijs}\bigg\rfloor \leq \ \sum_{s \in \sku}\xvar_{ijs} \leq \bigg\lceil \sum_{s \in \sku} X^{Rel}_{ijs}\bigg\rceil, \quad (i,j)\in \mov\\[0cm]
%
% UB-LB Rows1
%&\bigg\lfloor \sum_{\substack{s \in \sku,\,i:(i,j) \in \mov}} X^{Rel}_{ijs}\bigg\rfloor \leq \ \sum_{\substack{s \in \sku,\,i:(i,j) \in \mov}}\xvar_{ijs} \leq \bigg\lceil \sum_{\substack{s \in \sku,\,i:(i,j) \in \mov}} X^{Rel}_{ijs}\bigg\rceil, \quad j\in\fac\\[0cm]
%
% UB-LB Columns1
%&\bigg\lfloor \sum_{\substack{s \in \sku,\,j:(i,j) \in \mov}} X^{Rel}_{ijs}\bigg\rfloor \leq \ \sum_{\substack{s \in \sku,\,j:(i,j) \in \mov}}\xvar_{ijs} \leq \bigg\lceil \sum_{\substack{s \in \sku,\,j:(i,j) \in \mov}} X^{Rel}_{ijs}\bigg\rceil, \quad i\in\fac\\[0cm]
%
% UB-LB RefTotal
&\Big\lfloor B_{is}\Big\rfloor \leq \sum_{\substack{j:(j,i) \in \mov}}\xvar_{jis}- \sum_{\substack{j:(i,j) \in \mov}}\xvar_{ijs} \leq \Big\lceil B_{is}\Big\rceil, \quad i\in\fac \label{eq:balsku}\\[0cm]
%
% UB-LB Total
%&\bigg\lfloor \sum_{\substack{(i,j)\in\mov,\,s \in \sku}} X^{Rel}_{ijs}\bigg\rfloor \leq \ \sum_{\substack{(i,j)\in\mov,\,s \in \sku}}\xvar_{ijs} \leq \bigg\lceil \sum_{\substack{(i,j)\in\mov,\,s \in \sku}}X^{Rel}_{ijs}\bigg\rceil\\[0cm]
%
% Integrity
&\xvar_{ijs} \in \mathbb{Z^{+}}, \quad (i,j)\in\mov. \notag
\end{align}

As mentioned above, despite the integrality constraints in the $\xvar$ variables, solving the relaxed version of \rpm{s} already results in integer solutions. The objective function is discussed in Section~\ref{subsec:roundingcost} below. Constraints in~\eqref{eq:round_fac} ensure that the resulting values for the $\xvar_{ijs}$ variables are rounded counterparts of the $X^{Rel}_{ijs}$ parameters, defined from the solution of \rtp. Then, \eqref{eq:round_bal_fac} ensures that the total number of items sent from each facility~$i$ is a rounding of $Z^S_{is}$. Similarly, \eqref{eq:round_fac_bal} ensures that the total number of items received by each facility~$i$ is a rounding of $Z^R_{is}$. Constraints~\eqref{eq:round_fac}-\eqref{eq:round_fac_bal} are the classic matrix rounding constraints. In our setting, however, the additional set of constraints given by~\eqref{eq:balsku} is also needed to ensure that the net balance of SKU~$s$ at facility~$i$, that is, $\sum_{j:(j,i)\in\mov}\xvar_{jis}-\sum_{j:(i,j)\in\mov}\xvar_{ijs}$,  is a rounded counterpart of the $B_{is}$ parameter. Otherwise, this net balance may also be $\lfloor B_{is}\rfloor-1$ and $\lceil B_{is}\rceil+1$, which has two main problems:
\begin{itemize}
	\item Consider the following situation for a given SKU~$s$. An outlet has a fixed demand of 2 units and an initial stock of 1 unit. In the solution of the relaxed problem it sends $0.3$ units and receives $1.5$ units. This results in a final stock of $2.2$ units, which is enough to meet the fixed demand. Now, constraints~\eqref{eq:round_fac}-\eqref{eq:round_fac_bal} allow for a solution under which the outlet both sends and receives $1$ unit, resulting in an infeasible final stock of $1$ unit. The net balance in the relaxed version was $1.2$, but in the rounded solution this net balance becomes $0$.
	
	Although this particular infeasibility problem does not arise under the constraints in~\eqref{cons:T_dsiyn}, it would be a concern if they were replaced by the weaker version in~\eqref{cons:T_sendifav}. 	
	
	\item Equations \eqref{cons:T_obj2} and \eqref{cons:sfin} in problem \tp imply that the net balance of SKU $s$ in facility~$i$ is also important for the objective function. This is specially so for large values of $\alpha$, the penalty parameter for unmet variable demand. Constraints in~\eqref{eq:balsku} aim at producing solutions of the rounding problem whose objective function is closer to the objective function of the relaxed solution in problem \rtp.
\end{itemize}

\subsubsection{Objective function in problem \rp}
\label{subsec:roundingcost}

As discussed in Section~\ref{subsec:relaxed}, the main limitation of the \rtprp scheme is that, since the packing component of \tp problem is not present in \rp, extra packages may be required after the rounding takes place, with the corresponding increase in the cost of the resulting feasible solution of \tp. In order to mitigate this limitation, one can exploit the freedom in the choice of the $\hat{c}_{ijs}$ cost parameters that define the objective function of \rpm{s}, with $s\in \sku$. Since the rounding problems associated to the different SKUs are solved sequentially one can, when problem \rpm{s} is about to be solved, choose the $\hat{c}_{ijs}$ costs to incentivize the rounding up of transfers between facilities whose corresponding packages have more capacity left. In particular, we use costs of the form
\begin{equation}
\label{eq:rounding cost}
\hat{c}_{ijs}:=-\frac{\text{\footnotesize{``Total capacity of packages sent from $i$ to $j$''}}-\text{\footnotesize{``Total weight of items sent from $i$ to $j$''}}}{\text{\footnotesize{``Average package's cost from $i$ to $j$''}}}.
\end{equation}

Note that these costs won't be the same for all SKUs, since they must be updated after each \rpm{s} problem is solved.

Further, as we briefly discuss in Section~\ref{sec:TRP_num}, since the \rp solving stage is not demanding computationally, one can even solve it several times for different specifications of the $\hat{c}_{ijs}$ costs and return the best solution obtained in the process.

\subsubsection{Minimum cost flow formulation}

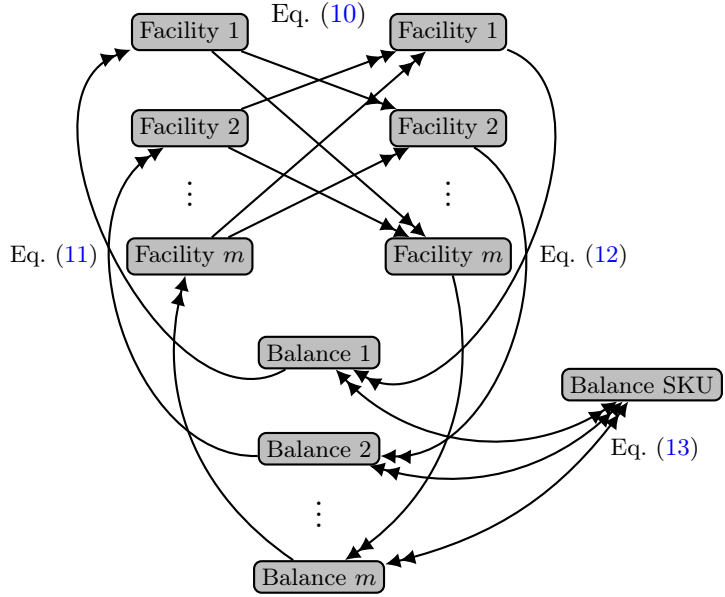
\begin{figure}[!htbp]
\centering
%\begin{subfigure}{.34\textwidth}
	%\centering
		%\begin{tikzpicture}[scale=0.85,>={Latex[length=2mm,width=2mm]},font=\footnotesize,every node/.style={font=\footnotesize},thick,estilo/.style={rectangle,draw,inner sep=1mm,rounded corners=1mm,fill=gray!50},theshift/.style={xshift=7.5cm,yshift=0cm}]
		%%%%%%Start grid
		%%\draw[step=.5cm,gray, thin, dashed] (-6,-4) grid (6,4);
		%%\foreach \x in {-6,-5,...,6}
		%%\draw (\x, -4) node[anchor=north, darkgray, fill=white] {{\tiny $\x$}};
		%%\foreach \y in {-4,-3,...,4}
		%%\draw (-6, \y) node[anchor=east, darkgray, fill=white] {{\tiny $\y$}};
		%%%%%%End grid
		%
		%%\node[rectangle,draw,inner sep=0.5mm,rounded corners=1mm,fill=gray!50] (fuente) at (0,0)  {ref1};
		%
		%\node[estilo] (fila1) at (3,1.5){Facility $1$};
		%\node[estilo] (fila2) at (3,0){Facility $2$};
		%\node[estilo] (fila3) at (3,-2){Facility $m$};
		%\node[estilo] (columna1) at (7,1.5){Facility $1$};
		%\node[estilo] (columna2) at (7,0){Facility $2$};
		%\node[estilo] (columna3) at (7,-2){Facility $m$};
		%
		%\node (dots1) at (3, -0.9){\Large{\vdots}};
		%\node (dots2) at (7, -0.9){\Large{\vdots}};		
%
		%\draw[->>] (fila1) -- (columna2);
		%\draw[->>] (fila1) -- (columna3);
		%\draw[->>] (fila2) -- (columna1);
		%\draw[->>] (fila2) -- (columna3);
		%\draw[->>] (fila3) -- (columna1);
		%\draw[->>] (fila3) -- (columna2);		
%
	%\end{tikzpicture}
	%\caption{Matrix rounding problem.}
	%\label{fig:matround}
%\end{subfigure}
%\begin{subfigure}{.6\textwidth}
	%\centering
		\begin{tikzpicture}[scale=0.85,>={Latex[length=2mm,width=2mm]},font=\footnotesize,every node/.style={font=\footnotesize},thick,estilo/.style={rectangle,draw,inner sep=1mm,rounded corners=1mm,fill=gray!50},theshift/.style={xshift=7.5cm,yshift=0cm}]
		%%%%%%%Start grid
		%\draw[step=.5cm,gray, thin, dashed] (-6,-4) grid (6,4);
		%\foreach \x in {-6,-5,...,6}
		%\draw (\x, -4) node[anchor=north, darkgray, fill=white] {{\tiny $\x$}};
		%\foreach \y in {-4,-3,...,4}
		%\draw (-6, \y) node[anchor=east, darkgray, fill=white] {{\tiny $\y$}};
		%%%%%%%End grid
		
		%\node[rectangle,draw,inner sep=0.5mm,rounded corners=1mm,fill=gray!50] (fuente) at (0,0)  {ref1};
		
		\node[estilo] (fila1) at (3,1.5){Facility $1$};
		\node[estilo] (fila2) at (3,0){Facility $2$};
		\node[estilo] (fila3) at (3,-2){Facility $m$};
		\node[estilo] (columna1) at (7,1.5){Facility $1$};
		\node[estilo] (columna2) at (7,0){Facility $2$};
		\node[estilo] (columna3) at (7,-2){Facility $m$};
		
		\node (dots1) at (3, -0.9){\Large{\vdots}};
		\node (dots2) at (7, -0.9){\Large{\vdots}};
		
		\node[estilo] (balance1) at (5,-3.5) {Balance 1} edge [bend left = 95,->>] (fila1) edge [bend right = 95,<<-] (columna1);
		\node[estilo] (balance2) at (5,-5){Balance 2} edge [bend left = 75,->>] (fila2) edge [bend right = 75,<<-] (columna2);
		\node[estilo] (balance3) at (5,-7){Balance $m$} edge [bend left = 35,->>] (fila3) edge [bend right = 35,<<-] (columna3);
		\node (dots2) at (5, -5.9){\Large{\vdots}};
		
		\node[estilo] (balanceref) at (10,-4){Balance SKU} edge [bend left = 40,<<->>] (balance1) edge [bend left = 30,<<->>] (balance2) edge [bend left = 20,<<->>] (balance3);
		
		\draw[->>] (fila1) -- (columna2);
		\draw[->>] (fila1) -- (columna3);
		\draw[->>] (fila2) -- (columna1);
		\draw[->>] (fila2) -- (columna3);
		\draw[->>] (fila3) -- (columna1);
		\draw[->>] (fila3) -- (columna2);	
		
		\draw[draw=none] (fila1) -- node[above=-0.1cm,sloped, font=\small]{Eq.~\eqref{eq:round_fac}} (columna1);
		\node (fac_bal)  at (9.1,-2) {Eq.~\eqref{eq:round_fac_bal}};
		\node (bal_fac)  at (0.9,-2) {Eq.~\eqref{eq:round_bal_fac}};
		\node (balsku)  at (10.2,-5) {Eq.~\eqref{eq:balsku}};

	\end{tikzpicture}
	%\caption{Modification of the matrix rounding problem.}
	%\label{fig:matroundmod}
%\end{subfigure}
\caption{Rounding scheme as a minimum cost flow problem.}
	\label{fig:rounding_scheme}
\end{figure}

We briefly outline how the rounding scheme in \rpm{s} can be formulated as a minimum cost flow problem, building upon the network flow representation of the matrix rounding problem discussed in Chapter~6.2 of \cite{Ahuja1993}. This formulation ensures the integrality of the solutions of the linear relaxation of \rpm{s}. The associated graph is represented in Figure~\ref{fig:rounding_scheme} which, with respect to the classic matrix rounding problem, adds all the balance nodes.

To complete the specification of the associated network flow problem we specify the lower and upper bounds associated to the flows through the different edges, along with the corresponding cost parameters. The flow of an edge between two facilities~$i$ and~$j$, $\xvar_{ijs}$, must lie in $\big[\lfloor X^{Rel}_{ijs}\rfloor, \lceil X^{Rel}_{ijs}\rceil \big]$ and each unit of flow has cost $\hat{c}_{ijs}$. For each facility~$i$, the flow of the edge from the balance node to the facility must lie in $\big [\lfloor Z^{S}_{is}\rfloor, \lceil Z^{S}_{is}\rceil \big ]$ and has zero cost; the flow balance constraint at this facility node requires that this flow is $\sum_{\substack{j:(i,j) \in \mov}}\xvar_{ijs}$. Similarly, for each facility~$i$, the flow of the edge from the facility to the balance node must lie in $\big[\lfloor Z^{R}_{is}\rfloor, \lceil Z^{R}_{is}\rceil \big]$ and has zero cost; the flow balance constraint at this facility node requires that this flow is $\sum_{\substack{j:(j,i) \in \mov}}\xvar_{jis}$. Finally, for each facility~$i$, the orientation of the flow in the edge joining $i$'s balance node with node balance SKU depends on the sign of $B_{is}$ and must lie in $\big[\lfloor \card{B_{is}}\rfloor, \lceil \card{B_{is}}\rceil \big]$, also with zero cost; the flow balance constraint at $i$'s balance node requires that this flow is $\sum_{\substack{j:(j,i) \in \mov}}\xvar_{jis}- \sum_{\substack{j:(i,j) \in \mov}}\xvar_{ijs}$.\footnote{The flow balance constraint at node ``Balance SKU'' is trivially satisfied, since it reduces to require that the overall net balance is zero: $\sum_{i \in \fac} \Big(\sum_{j:(j,i)\in\mov}\xvar_{jis}-\sum_{j:(i,j)\in\mov}\xvar_{ijs}\Big)=0$.}

\subsection{\TRPA. Phase 3: Packing}\label{subsec:phase3}

The final phase of the \TRPA, which is also needed by the direct solution via the \TPA approach, is to solve the packing problems \ppm{ij}, one for each $(i,j) \in \mov$. These problems were defined in Section~\ref{subsec:packingdef}, where its role to enhance the feasible solutions of \tp was also discussed.

Importantly, when seen as part of \TRPA, this additional phase may even help in reducing the cost of the solutions obtained after \rtprp has been run. The reason is that, in the rounding phase, when a new package is required between two facilities, we just ``myopically'' add the cheapest available package between them. When the corresponding packing problem is solved, the final solution contains the optimal package configuration given the items sent between those two facilities, regardless of the specific packages that might have been chosen in the previous phases. Importantly, the solution of the packing phase also provides, for each pair of facilities, the precise composition of each package to be sent between them.

\section{\TRPA: Numeric Analysis}
\label{sec:TRP_num}

In this section we present a thorough comparison of the performance of the solution approach \TRPA with the more direct \TPA. The analysis is specially focused on the scalability of the solution times for large instances, those for which the \TRPA was conceived. Before discussing the numeric results, we present a brief discussion on the specifics of how \TPA and \TRPA approaches have been implemented and how the different test sets have been generated.

\subsection{Framework for the comparison between \TPA and \TRPA}\label{subsec:frame}

Both \TPA and \TRPA solution methods have been implemented in Python 3.8 \citep{python3} and all the MILP problems have been solved with Gurobi 8.0.0 \citep{GurobiOptimization2020},  interfaced via its Python API. In all instances, the maximum solution time for both \TPA and \TRPA has been set to 5 minutes. In addition, the largest instances, discussed in Section~\ref{sec:huge}, have also been solved with a 30-minute time limit. All the executions reported in this section have been performed on the supercomputer Finisterrae~II, located at Galicia Supercomputing Centre (CESGA).\footnote{For further information refer to \url{https://www.cesga.es/}.} Specifically, we used computational nodes powered with 2 deca-core Intel Haswell 2680v3 CPUs with 128GB of RAM connected through an Infiniband FDR network, and 1TB of hard drive.

Next, we briefly discuss how the rounding and packing schemes have been implemented: 
\begin{description}
\item[Relaxed transferring, \rtp.]  \TRPA is run with four values of $\delta$: $0.85$, $0.9$, $0.95$, and~$1$.%\footnote{Values of $\delta$ below $0.8$ were not competitive in our test sets.} %Each of the resulting solutions is then compared with the solution of \TPA.
\item[Rounding, \rp.] As discussed in Section~\ref{subsec:roundingcost}, when defining the problems \rpm{s} there is complete freedom to define the $\hat{c}_{ijs}$ cost parameters of the rounding scheme. In each call to \TRPA, the full \rp problem is solved repeatedly for a maximum of 50 runs,\footnote{Recall that solution times in \rp are almost negligible when compared to the solution times of the other phases of the algorithm.} returning the best solution found, according to the following rules: i)~the cost parameters for the first run are of the form of Equation~\eqref{eq:rounding cost}, with these parameters being randomly perturbed in subsequent runs, ii)~the SKUs $s\in \sku$ are decreasingly sorted according to their weights in the first run and randomly sorted in subsequent runs, iii)~if after solving a certain  \rpm{s} problem an additional package is needed between a pair of facilities, one unit of the cheapest available package is added, iv)~the procedure stops if a solution is found in which no additional package is required and, moreover, its cost in \tp is within a given tolerance of the cost of the solution of problem \rtp,\footnote{In particular, all the impact of the rounding in the cost function must come from the variable demand term in Equation~\eqref{cons:T_obj2}.} and v)~it also stops if solutions with the same (best) associated cost in \tp are found in 5 consecutive runs. 
\item[Packing, \pp.] The availabilities for the different packages in $\pack$, $\avap_p$, are set to be large enough to ensure that there is excess capacity in all types of packages.
\end{description}

Given that the packing problem is common to both \TRPA and \TPA approaches, when comparing the computational times of both solution methods we just focus on the running times of \rtprp versus those of \tp.

\subsection{Generation of the test sets}\label{sec:TS}

For the generation of the test sets we rely again on the simulator briefly discussed in Section~\ref{sec:numeric_policies} and further explained in Appendix~\ref{sec:appendix_generator}. We have generated four different test sets, depending on the size of the resulting problems:
\begin{description}
\item[Small-TS.] 100 instances. We generate 25 instances with $\card{\sku}=10$, $\card{\pack}=2$, $\card{\ware}=1$, $\card{\outl}=10$, and a total stock in the network of $1000$ items. From each of these baseline instances we generate 4 instances by giving $\alpha$ parameter the following values: $0$, $0.1$, $10$, and $1000$. Each of these instances has $1430$ variables and $430$ constraints; \tp has $1320$ integer variables whereas \rtp has just $110$.
\item[Medium-TS.] 100 instances. We generate 25 instances with $\card{\sku}=30$, $\card{\pack}=4$, $\card{\ware}=1$, $\card{\outl}=30$, and a total stock in the network of $9000$ items. Again, four variations of each were created by giving $\alpha$ parameter the values $0$, $0.1$, $10$, and $1000$. Each of these instances has $32550$ variables and $3690$ constraints; \tp has $31620$ integer variables whereas \rtp has just $3720$.
\item[Large-TS.] 100 instances. We generate 25 instances with $\card{\sku}=100$, $\card{\pack}=4$, $\card{\ware}=1$, $\card{\outl}=100$, and a total stock in the network of $100000$ items. Again, four variations of each were created by giving $\alpha$ parameter the values $0$, $0.1$, $10$, and $1000$. Each of these instances has $1060500$ variables and $40300$ constraints; \tp has $1050400$ integer variables whereas \rtp has just $40400$.
\item[Extra-large-TS.] Finally, we generate a last test set with instances of different sizes to push both \TRPA and \TPA to their limits, and whose results are discussed in Section~\ref{sec:huge}. We generate 10 groups of 25 instances in each of them. Each of these 10 groups was generated as the above ones, with $\card{\ware}=1$, parameter $\alpha$ taking just value~0, and the specifications for $\card{\sku}$, $\card{\pack}$, and $\card{\outl}$ given in Table~\ref{tab:huge_cases}, which also contains information about the sizes of the resulting instances. 
\end{description}

\begin{table}[!htbp]
    \centering
    \begin{tabular}{lllllllll}
        \multirow{2}{*}{Group} & \multirow{2}{*}{$\card{\sku}$} & \multirow{2}{*}{$\card{\pack}$} & \multirow{2}{*}{$\card{\outl}$} & Overall   & \multirow{2}{*}{Constraints}  &\multirow{2}{*}{Variables} & \tp integer  & \rtp  integer   \\ 
				        &  &  &  & stock    &  & & variables & variables  \\ 
        \hline
        G1     & 80         & 2        & 80    & 64000   & 25840 & 537840 & 531360 & 12960   \\
        G2     & 100        & 2        & 100   & 100000  & 40300 & 1040300 & 1030200 & 20200\\
        G3     & 120        & 2        & 120   & 144000  & 57960 & 1785960 & 1771440 & 29040\\
        G4     & 140        & 2        & 140   & 196000  & 78820 & 2822820 & 2803080 & 39480 \\
        G5     & 160        & 2        & 160   & 256000  & 102880  & 4198880 & 4173120 & 51520 \\
        G6     & 180        & 2        & 180   & 324000  & 130140  & 5962140 & 5929560 & 65160\\
        G7     & 200        & 2        & 200   & 400000  & 160600  & 8160600 & 8120400 & 80400\\
        G8     & 220        & 2        & 220   & 484000  & 194260  & 10842260 & 10793640 & 97240\\
        G9     & 240        & 2        & 240   & 576000  & 231120  & 14055120 & 13997280 & 115680 \\
        G10    & 260        & 2        & 260   & 676000  & 271180 & 17847180 & 17779320  & 135720\\
    \end{tabular}
    \caption{Description of the 10 groups of instances in Extra-large T-S.}
    \label{tab:huge_cases}
\end{table}

\subsection{Results for Small-TS, Medium-TS, and Large-TS}
\label{sec:results_instances_sml}

\begin{figure}[!htbp]
\begin{subfigure}{.46\textwidth}
  \centering
	\vspace{0.05cm}
  \includegraphics[width=\textwidth,viewport=15mm 18mm 219mm 123mm,clip]{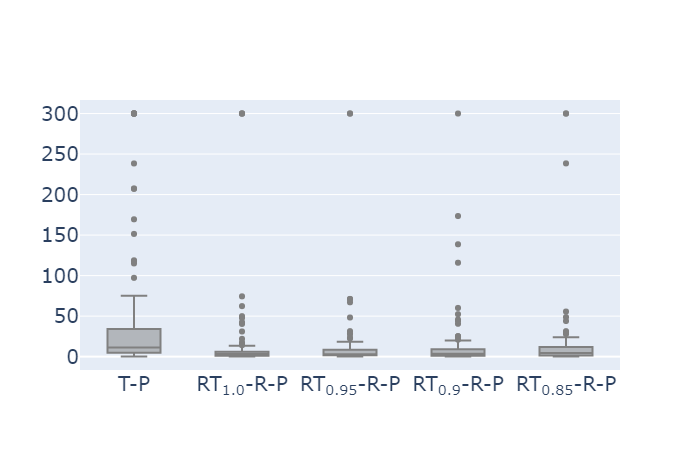}
	%\vspace{0cm}
  \caption{Box plots (in seconds).}
  \label{fig:time_small_trans}
\end{subfigure}
\hspace{0.2cm}
\begin{subfigure}{.54\textwidth}
  \centering
	\includegraphics[width=\textwidth,viewport=16mm 20mm 241mm 130mm,clip]{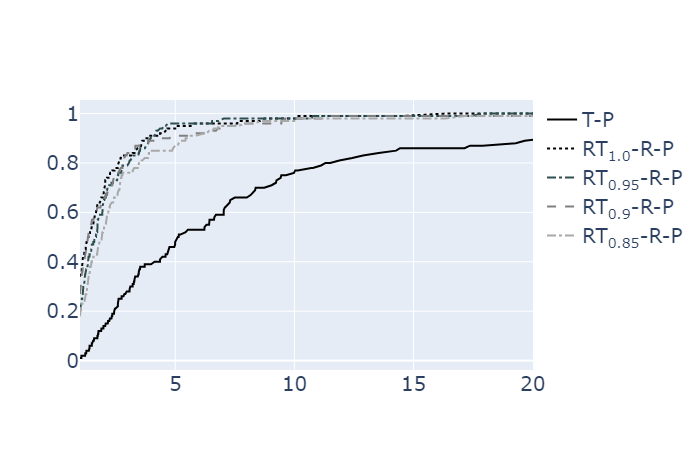}
	\caption{Performance profile.}
	\label{fig:Small_totalTimeDS_PerformanceProfile}
\end{subfigure}
%
%\begin{subfigure}{\textwidth}
  %\centering
   %include second image
  %\includegraphics[width=.48\textwidth,viewport=19mm 22mm 219mm 130mm,clip]{images/Large_totalTimeTransDSBoxPlot.png}  
  %\caption{Large instances}
  %\label{fig:time_large_trans}
%\end{subfigure}
\caption{Computational times in \tp and \rtp problems in Small-TS.}
\label{fig:time_boxplot_Trans}
\end{figure}

We start discussing the computational times. Figure~\ref{fig:time_small_trans} represents the running times in Small-TS associated with the solution by Gurobi of problems \tp and \rtp, with $\delta \in \{1,0.95,0.9,0.85\}$. As expected, \rtp is notably faster, although we can see that, for all values of $\delta$, there are instances that were not solved to optimality after 5~minutes (300~seconds) even in this set of relatively small instances. In Figure~\ref{fig:Small_totalTimeDS_PerformanceProfile} we represent the corresponding performance profiles \citep{Dolan2002}, which reinforce the results in Figure~\ref{fig:time_small_trans} and provide some additional information. For instance, we can see that the probability that \tp solves a given instance more than 10 times slower that the best \rtp configuration is larger than $0.2$. On the other hand, this probability is very close to zero for all \rtp configurations. We do not represent the corresponding plots for  Medium-TS and Large-TS, since in both of them almost all instances reached the 5-minute time limit without closing the optimality gap for both \tp and all four \rtp configurations. 

\begin{figure}[!htbp]
\begin{subfigure}{.32\textwidth}
  \centering
  \includegraphics[width=\textwidth,viewport=5mm 12mm 258mm 177mm,clip]{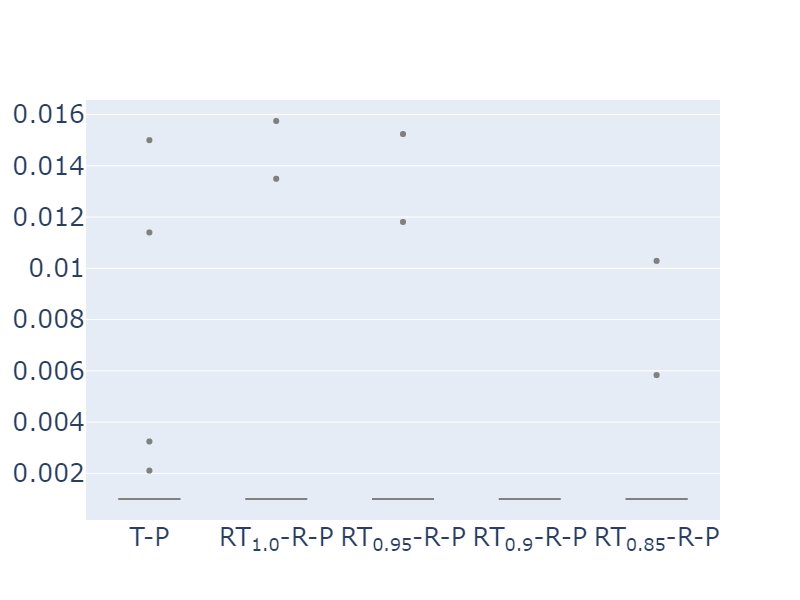}
  \caption{Small-TS.}
  \label{fig:gap_small_trans}
\end{subfigure}
\hspace{0.05cm}
\begin{subfigure}{.32\textwidth}
  \centering
  \includegraphics[width=\textwidth,viewport=5mm 12mm 258mm 177mm,clip]{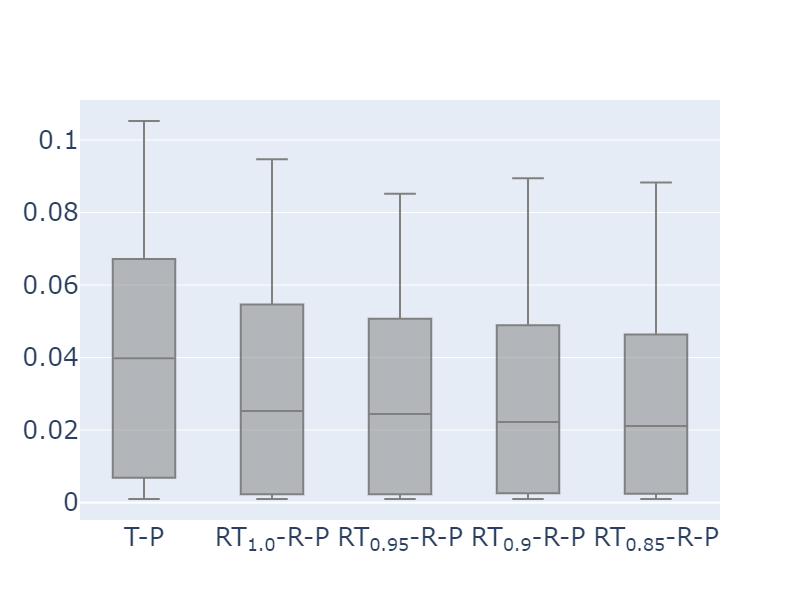}  
  \caption{Medium-TS.}
  \label{fig:gap_medium_trans}
\end{subfigure}
\hspace{0.05cm}
\begin{subfigure}{.32\textwidth}
  \centering
  % include second image
  \includegraphics[width=\textwidth,viewport=5mm 12mm 258mm 177mm,clip]{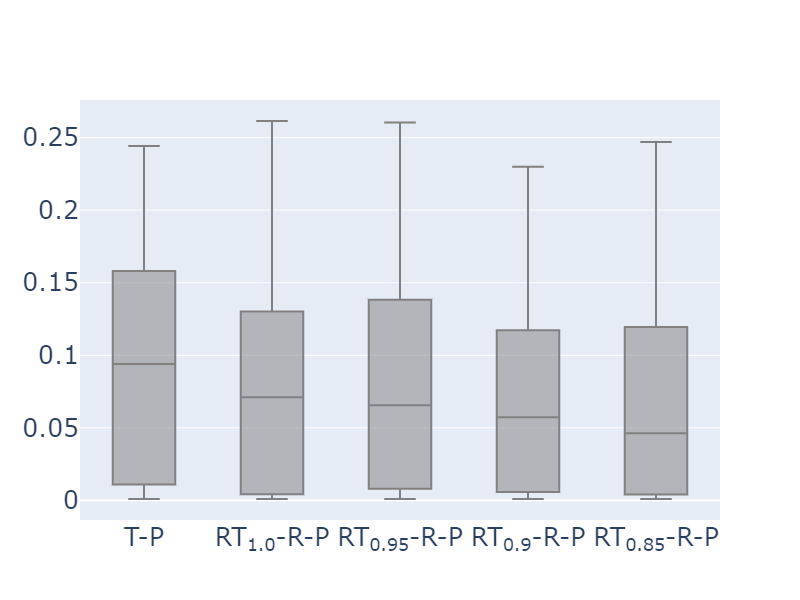}  
  \caption{Large-TS.}
  \label{fig:gap_large_trans}
\end{subfigure}
\caption{Average relative optimality gaps after \tp and \rtp, with $\delta \in \{1,0.95,0.9,0.85\}$.}
\label{fig:relgap_boxplot_Trans}
\end{figure}

We now move to a comparison of the average relative optimality gaps obtained by Gurobi after 5~minutes in \tp and \rtp. Note that here we expect to see that \rtp is closer to convergence, since \tp has many more integer variables and closing the gap should therefore be harder. This can be seen in the box plots in Figure~\ref{fig:relgap_boxplot_Trans}. Box plots for Small-TS just show that almost all instances are solved to optimality in this test set. Average gaps increase for Medium-TS and increase even further for Large-TS. These gaps are noticeably larger for \tp than for the different \rtp configurations, which have no such a clear pattern between them.

So far we have just seen that the results from these test sets confirm that the relaxed versions of \tp are easier to solve than \tp. Yet, this is not enough to justify that \TRPA can be competitive with \TPA. There are two main things that need to be discussed. First, the computation times of the rounding phase need to be small with respect to the solution times of \tp and \rtp, so that the computational overhead is not significant. Second, the quality of the final solutions obtained after applying the \rtprppp scheme must be comparable to those obtained with the direct solution via \tppp.

\begin{figure}[!htbp]
\begin{subfigure}{.32\textwidth}
  \centering
  \includegraphics[width=\textwidth,viewport=5mm 15mm 255mm 180mm,clip]{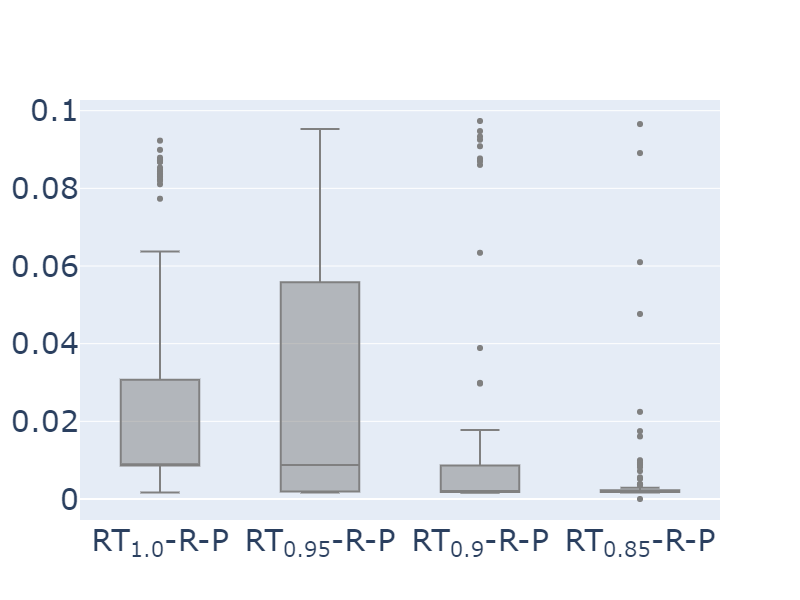}
  \caption{Small-TS.}
  \label{fig:time_small_round}
\end{subfigure}
\hspace{0.05cm}
\begin{subfigure}{.32\textwidth}
  \centering
  \includegraphics[width=\textwidth,viewport=5mm 15mm 255mm 180mm,clip]{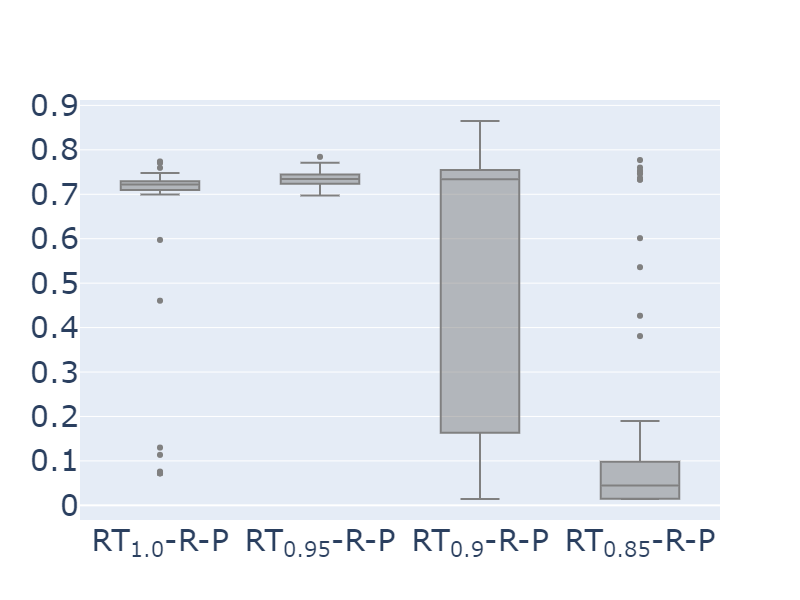}  
  \caption{Medium-TS.}
  \label{fig:time_medium_round}
\end{subfigure}
\hspace{0.05cm}
\begin{subfigure}{.32\textwidth}
  \centering
  % include second image
  \includegraphics[width=\textwidth,viewport=12mm 15mm 255mm 180mm,clip]{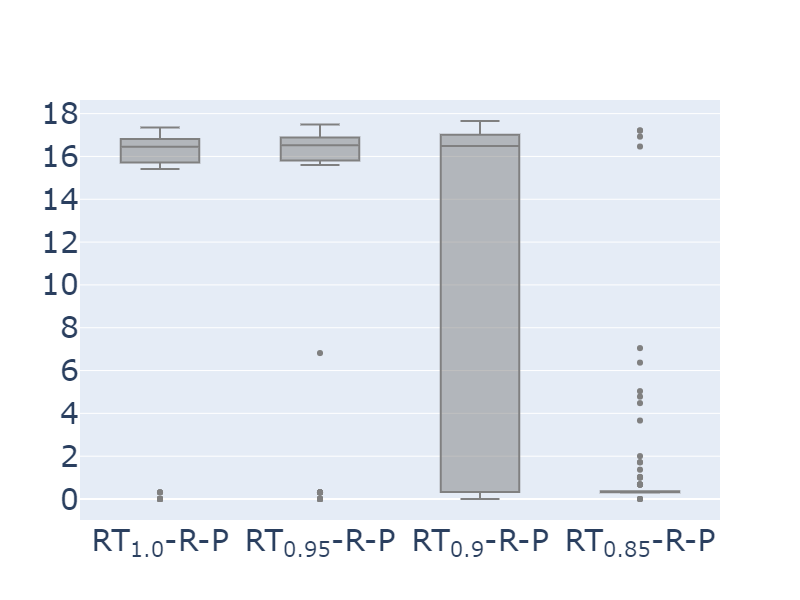}  
  \caption{Large-TS.}
  \label{fig:time_large_round}
\end{subfigure}
\caption{Computational times for rounding phase (in seconds).}
\label{fig:time_boxplot_Round}
\end{figure}

Figure~\ref{fig:time_boxplot_Round} contains box plots of the running times of the rounding phase, problem \rp, for the three test sets under analysis. We can see that this phase is indeed not demanding computationally since, in the worst cases, running times for Small-TS, Medium-TS, and Large-TS are under 0.1, 1, and 20 seconds, respectively. Interestingly, it seems that rounding problems are easier to solve for smaller values of $\delta$. This is natural, since the smaller $\delta$ is the easier it should be to find a rounding that does not require extra packages, prompting an early termination of the rounding loop described in Section~\ref{subsec:frame}. Note that this does not mean that smaller values of $\delta$ should be preferable, since they will require more packages in phase \rtp. 

\begin{figure}[!htbp]
\centering
\begin{subfigure}{.55\textwidth}
  \centering
  \includegraphics[width=\textwidth,viewport=15mm 20mm 280mm 180mm,clip=true]{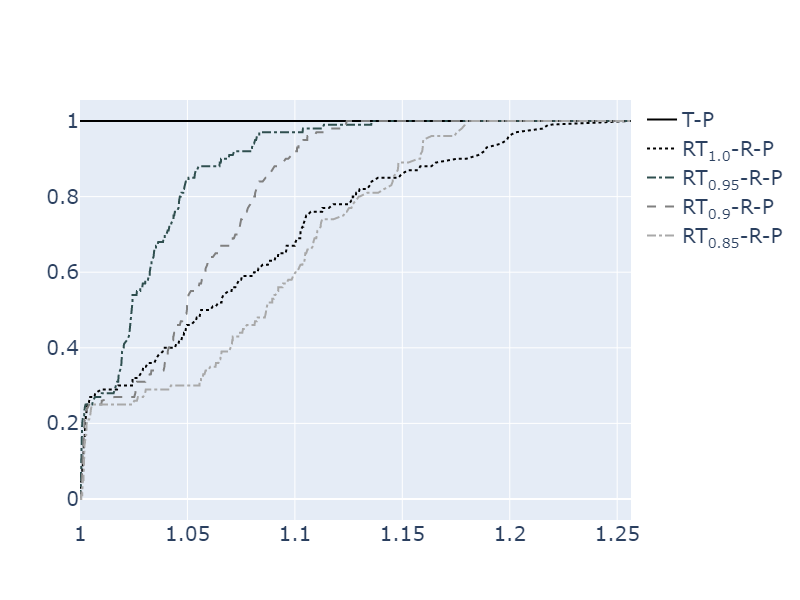}
  \caption{Small-TS after Rounding (\rp). \hspace*{7em}}
  \label{fig:obj_small_round}
\end{subfigure}
\begin{subfigure}{.43\textwidth}
  \centering
  \includegraphics[width=\textwidth,viewport=15mm 20mm 220mm 180mm,clip=true]{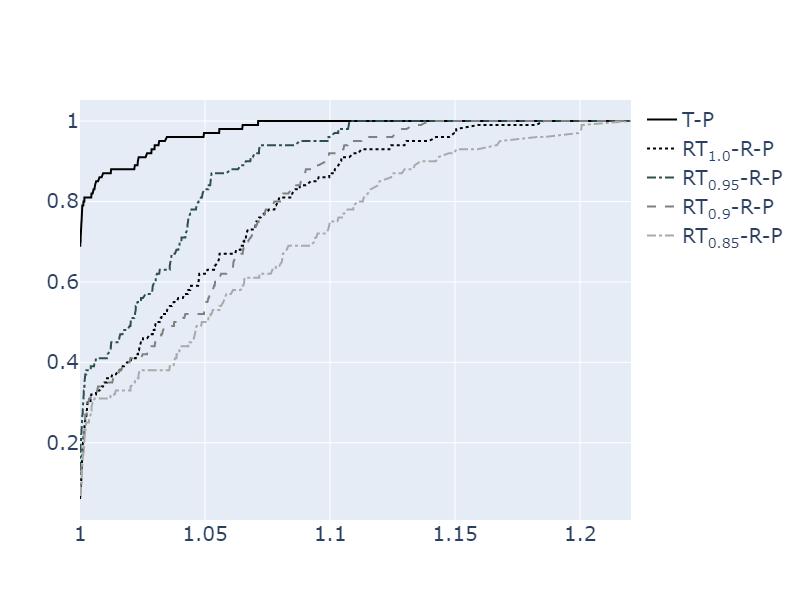}  
  \caption{Small-TS after Packing (\pp).}
  \label{fig:obj_small_pack}
\end{subfigure}
\begin{subfigure}{.55\textwidth}
  \centering
  \includegraphics[width=\textwidth,viewport=15mm 20mm 280mm 180mm,clip=true]{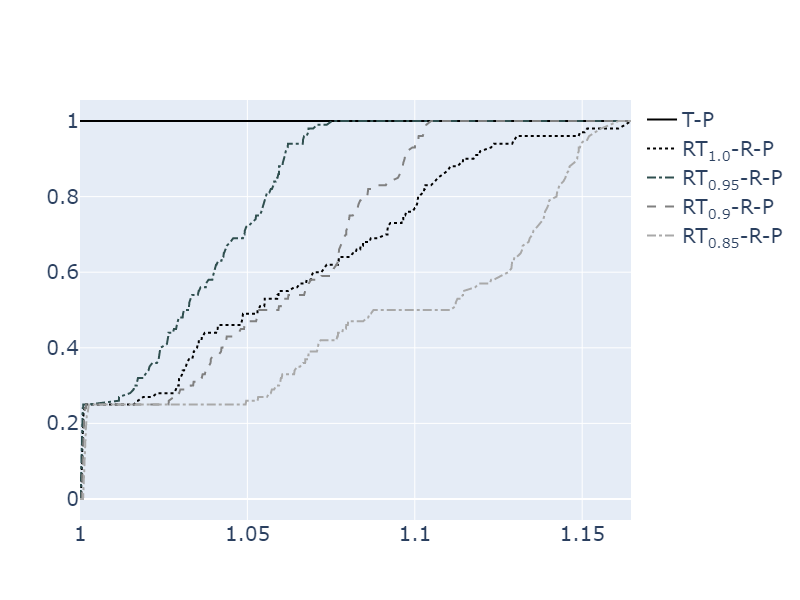}
  \caption{Medium-TS after Rounding (\rp). \hspace*{7em}}
  \label{fig:obj_medium_round}
\end{subfigure}
\begin{subfigure}{.43\textwidth}
  \centering
  \includegraphics[width=\textwidth,viewport=15mm 20mm 220mm 180mm,clip=true]{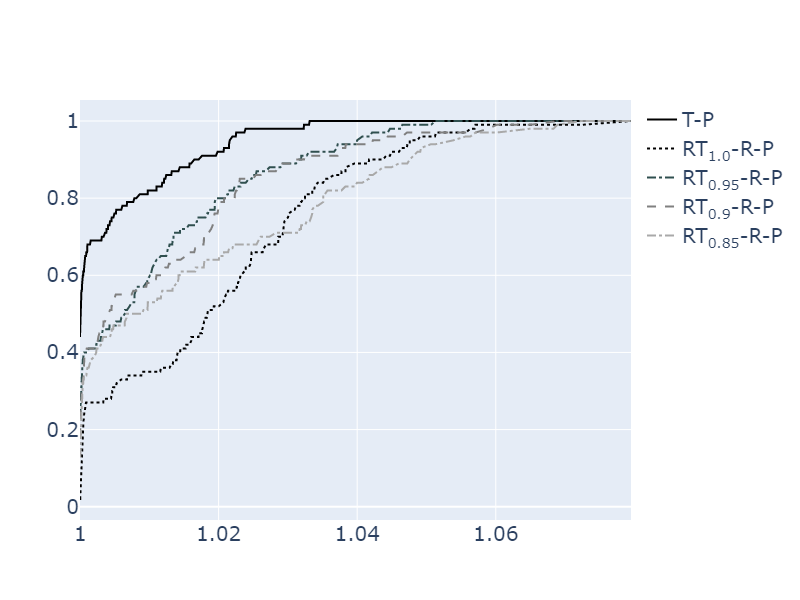}  
  \caption{Medium-TS after Packing (\pp).}
  \label{fig:obj_medium_pack}
\end{subfigure}
\begin{subfigure}{.55\textwidth}
  \centering
  \includegraphics[width=\textwidth,viewport=15mm 20mm 280mm 180mm,clip=true]{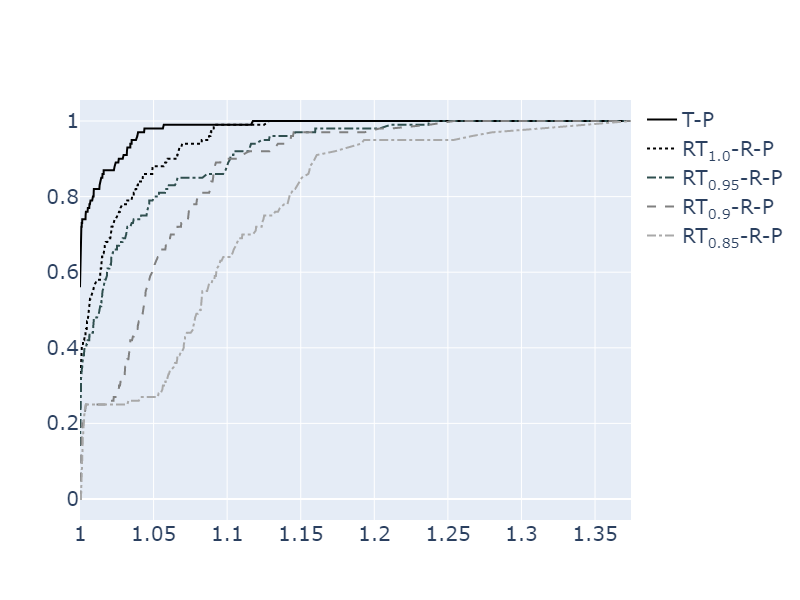}
  \caption{Large-TS after Rounding (\rp). \hspace*{7em}}
  \label{fig:obj_large_round}
\end{subfigure}
\begin{subfigure}{.43\textwidth}
  \centering
  \includegraphics[width=\textwidth,viewport=15mm 20mm 220mm 180mm,clip=true]{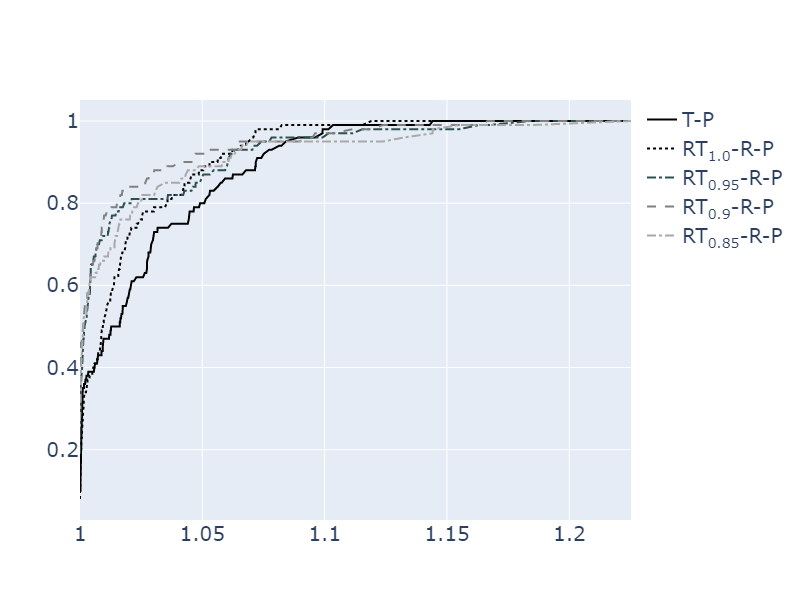}  
  \caption{Large-TS after Packing (\pp).}
  \label{fig:obj_large_pack}
\end{subfigure}
\caption{Performance profiles for the objective function for the different test sets.}
\label{fig:obj_perf_round_pack}
\end{figure}

Next, we present the most important part of the analysis, the comparison of the objective function of the solutions obtained by \TPA and \TRPA. The corresponding performance profiles are represented in Figure~\ref{fig:obj_perf_round_pack}. The graphs on the left correspond with the objective function of the solutions obtained after solving \tp and \rtprp, respectively, for \TPA and \TRPA algorithms. As expected, the solutions obtained by \tp are better than those obtained by \rtprp, regardless of the value of $\delta$. Indeed, in all instances of Small-TS and Medium-TS the direct solution of \tp delivers the best solution. In Large-TS \tp is still better, but it delivers the best solution in less than 60\% of the instances. This would not be possible if the instances were solved to optimality but, as we showed in Figure~\ref{fig:relgap_boxplot_Trans}, the optimality gaps for \tp in Large-TS are normally above 10\%, which makes it possible for the \rtprp scheme to find better solutions in some instances. 

% Funcion para ajustar coordenadas del viewport
% https://tex.stackexchange.com/questions/111401/draw-a-grid-for-obtaining-coordinates-to-pass-to-includegraphics-viewport-optio
%\newpage
%\begin{overpic}[scale=0.5,unit=1mm,grid,tics=5]{images/Small_objPackDS_PerformanceProfile.png}
  %\put(50,50){\makebox[0pt]{\colorbox{blue}{\color{white}\Huge Tiger}}}
	%\put(0,0){\makebox[0pt]{\colorbox{blue}{X}}}
%\end{overpic}
%\fbox{\includegraphics[viewport=25mm 20mm 170mm 130mm,scale=0.5,clip]{images/Small_objRoundDS_PerformanceProfile.png}}

%\clearpage

Despite the above observations, the real test to compare \TPA and \TRPA must come from the solutions after the full \tppp and \rtprppp schemes have been run. In particular, as we discussed in Section~\ref{subsec:phase3}, for \TRPA the packing phase may even help to ``correct'' some myopic package choices made during the \rtprp phases. This may (partially) compensate for the extra packages that may be required to meet the packing constraints in the \ppm{ij} problems. Figures~\ref{fig:obj_small_pack} and~\ref{fig:obj_medium_pack} show that, for Small-TS and Medium-TS, \TPA is better than all configurations of \TRPA. Yet, the margin is narrower than it was in Figures~\ref{fig:obj_small_round} and~\ref{fig:obj_medium_round}: after running \pp, the best \TRPA configurations find the best solution in around 40\% of the instances, whereas after running \rp this number was around 20\%. More importantly, Figure~\ref{fig:obj_large_pack} shows that for the instances in Large-TS the situation completely reverses and all \TRPA configurations seem to outperform the direct solution via \TPA.

Given the above discussion one may conclude that, once we have to settle for approximate solutions in the transferring phase, the \TRPA scheme may be preferable to \TPA. To provide some additional support to this claim and given that the performance of \TRPA relative to \TPA seems to improve as problem size increases, in Section~\ref{sec:huge} we have made one last computational study using the instances in Extra-large-TS. In this study we show that, given some pre-specified time limits, for very large problems \TRPA will be able to obtain feasible solutions more consistently than \TPA.

Yet, before moving to Section~\ref{sec:huge}, we want to briefly discuss the behavior of \TRPA as a function of parameter $\delta$. Looking again at the left plots in Figure~\ref{fig:obj_perf_round_pack}, maybe the only clear pattern after the rounding phase is that $\delta=0.85$ is outperformed by the other choices for $\delta$. One may even argue that $\delta=0.95$ might be the best choice, at least for Small-TS and Medium-TS. Yet, as we argued above, the real test is to make this comparison on the right plots, once the full \TRPA has been run. Interestingly, although $\delta=0.95$ is still arguably a good choice, the situation is not so clear anymore, with $\delta=0.9$ exhibiting a very good performance for both Medium-TS and Large-TS. Interestingly, for these two test sets it is also worth noting the improvement of the configuration with $\delta=0.85$ once the packing phase has been run. 

In order to get a better understanding of the trade-offs associated with the choice of parameter~$\delta$ in the \TRPA algorithm, we present in Table~\ref{tab:packages} a more detailed analysis of the evolution of both \TPA and \TRPA in their different phases. We focus on the total number of packages sent and the part of the objective function associated with the transportation costs, given by term \eqref{cons:T_obj1} in \tp. More precisely, Table~\ref{tab:packages} contains, for each test set, the following information:
\begin{description}
	\item[Packages sent.] For each phase of the algorithms and each of the 5 solving schemes under analysis, we compute the average of the total number of packages sent. For \TPA there is no change after \rp, since this phase is not present in this algorithm. In parentheses we present the increase with respect to the previous phase.
	\item[Transportation costs.] For each phase of the algorithms and each of the 5 solving schemes, we compute the average transportation costs, given by $\sum_{(i,j)\in\mov, \, p\in\pack} \cost_{ijp} \yvar_{ijp}$. The table contains, for each phase, the result of dividing each of these costs by the minimum cost among the 5 solving schemes. In parentheses we present the increase of the underlying transportation cost with respect to the previous phase.
\end{description}

\begin{table}[!htbp]
{\small
\begin{tabular}{|c|lllll|}
\hline \hline
  \multicolumn{6}{|c|}{\rule{0pt}{13pt} {\normalsize \textbf{Packages and transportation costs in Small-TS}}}                                                                     \\ \hline \hline
 \textbf{Packages}                                         &   \TPA              & \TRPAdelta{1}              & \TRPAdelta{0.95}          & \TRPAdelta{0.9}           & \TRPAdelta{0.85}          \\ \hline
\textbf{\tp/\rtp}      & 31.34             & 31.22             & 32.34            & 33.65            & 35.14            \\ 
\textbf{\rp}          & 31.34 {\footnotesize(0\%)}       & 35.79 {\footnotesize(14.64\%)}   & 32.92 {\footnotesize(1.79\%)}   & 33.92 {\footnotesize(0.8\%)}    & 35.25 {\footnotesize(0.31\%)}   \\ 
\textbf{\pp}           & 34.95 {\footnotesize(11.52\%)}   & 37.18 {\footnotesize(3.88\%)}    & 36.16 {\footnotesize(9.84\%)}   & 37.14 {\footnotesize(9.49\%)}   & 38.16 {\footnotesize(8.26\%)}   \\ \hline
 \textbf{Tr. cost}                                         &   \TPA              & \TRPAdelta{1}              & \TRPAdelta{0.95}          & \TRPAdelta{0.9}           & \TRPAdelta{0.85}          \\ \hline 
\textbf{\tp/\rtp} & 1.0069            & 1                 & 1.0401           & 1.0812           & 1.1304           \\ 
\textbf{\rp}     & 1 {\footnotesize(0\%)}           & 1.1139 {\footnotesize(12.16\%)}  & 1.0466 {\footnotesize(1.32\%)}  & 1.0795 {\footnotesize(0.53\%)}  & 1.1248 {\footnotesize(0.19\%)}  \\ 
\textbf{\pp}      & 1 {\footnotesize(9.96\%)}        & 1.0515 {\footnotesize(3.8\%)}    & 1.026 {\footnotesize(7.8\%)}    & 1.0466 {\footnotesize(6.61\%)}  & 1.0678 {\footnotesize(4.39\%)}  \\ \hline \hline
\multicolumn{6}{|c|}{\rule{0pt}{13pt} {\normalsize \textbf{Packages and transportation costs in Medium-TS}}}                                                                    \\ \hline \hline
 \textbf{Packages}                                         &   \TPA              & \TRPAdelta{1}              & \TRPAdelta{0.95}          & \TRPAdelta{0.9}           & \TRPAdelta{0.85}          \\ \hline
\textbf{\tp/\rtp}      & 237.42            & 228.76            & 239.41           & 251.36           & 264.31           \\ 
\textbf{\rp}          & 237.42 {\footnotesize(0\%)}      & 259.88 {\footnotesize(13.6\%)}   & 246.78 {\footnotesize(3.08\%)}  & 252.59 {\footnotesize(0.49\%)}  & 264.45 {\footnotesize(0.05\%)}  \\ 
\textbf{\pp}           & 294.21 {\footnotesize(23.92\%)}  & 299.66 {\footnotesize(15.31\%)}  & 294.9 {\footnotesize(19.5\%)}   & 294.21 {\footnotesize(16.48\%)} & 296.47 {\footnotesize(12.11\%)} \\ \hline
 \textbf{Tr. cost}                                         &   \TPA              & \TRPAdelta{1}              & \TRPAdelta{0.95}          & \TRPAdelta{0.9}           & \TRPAdelta{0.85}          \\ \hline
\textbf{\tp/\rtp} & 1.02              & 1                 & 1.0479           & 1.0996           & 1.1574           \\ 
\textbf{\rp}     & 1 {\footnotesize(0\%)}           & 1.0765 {\footnotesize(9.8\%)}    & 1.0474 {\footnotesize(1.95\%)}  & 1.0814 {\footnotesize(0.31\%)}  & 1.1351 {\footnotesize(0.03\%)}  \\ 
\textbf{\pp}      & 1 {\footnotesize(18.17\%)}       & 1.0157 {\footnotesize(11.5\%)}   & 1.0035 {\footnotesize(13.22\%)} & 1.001 {\footnotesize(9.38\%)}   & 1.0042 {\footnotesize(4.54\%)}  \\ \hline \hline
\multicolumn{6}{|c|}{\rule{0pt}{13pt} {\normalsize \textbf{Packages and transportation costs in Large-TS}}}                                                                     \\ \hline \hline
 \textbf{Packages}                                         &   \TPA              & \TRPAdelta{1}              & \TRPAdelta{0.95}          & \TRPAdelta{0.9}           & \TRPAdelta{0.85}          \\ \hline
\textbf{\tp/\rtp}      & 2696.07           & 2594.66           & 2722.27          & 2810.47          & 2936.81          \\ 
\textbf{\rp}          & 2696.07 {\footnotesize(0\%)}     & 2718.75 {\footnotesize(4.78\%)}  & 2750.89 {\footnotesize(1.05\%)} & 2817.77 {\footnotesize(0.26\%)} & 2938.91 {\footnotesize(0.07\%)} \\ 
\textbf{\pp}           & 3047.36 {\footnotesize(13.03\%)} & 3009.06 {\footnotesize(10.68\%)} & 2993.27 {\footnotesize(8.81\%)} & 2970.61 {\footnotesize(5.42\%)} & 2979.6 {\footnotesize(1.38\%)}  \\ \hline
 \textbf{Tr. cost}                                         &   \TPA              & \TRPAdelta{1}              & \TRPAdelta{0.95}          & \TRPAdelta{0.9}           & \TRPAdelta{0.85}          \\ \hline
\textbf{\tp/\rtp} & 1.0156            & 1                 & 1.0499           & 1.0903           & 1.1431           \\ 
\textbf{\rp}     & 1 {\footnotesize(0\%)}           & 1.0167 {\footnotesize(3.25\%)}   & 1.0407 {\footnotesize(0.67\%)}  & 1.0753 {\footnotesize(0.17\%)}  & 1.1261 {\footnotesize(0.05\%)}  \\ 
\textbf{\pp}      & 1.0184 {\footnotesize(10.1\%)}   & 1.0115 {\footnotesize(7.55\%)}   & 1.0068 {\footnotesize(4.58\%)}  & 1.0001 {\footnotesize(0.55\%)}  & 1 {\footnotesize(-4\%)}         \\ \hline
\end{tabular}
}
\caption{Evolution of the packages required by the different solution schemes.}
\label{tab:packages}
\end{table}

As expected, \TRPAdelta{1}, delivers the best results after the transferring phase in the three test sets: minimum number of packages sent and lowest transportation costs, with both values increasing as $\delta$ decreases. %Recall that \TRPAdelta{1} dominates the solution provided by \TPA after the transferring phase, given that \rtpdelta{1} is just a relaxation of the integrality of the $\xvar_{ijs}$ variables in \tp. 

Next, we have that the number of extra packages required in the rounding phase is increasing in $\delta$, which is natural, given that reducing $\delta$ lead to solutions with more slack in the packages. Thus, the transportation costs of the different \TRPA configurations get closer to each other. Importantly, these costs also show, as we had already seen in the left plots of Figure~\ref{fig:obj_perf_round_pack}, that the solution provided by \TPA after \tp dominates all \rtprp solutions.

Finally, we analyze the impact of the packing phase, in which the \ppm{ij} problems are solved. This phase may call for additional packages for both \TPA and all \TRPA configurations. This increase is noticeably larger for \TPA in all test sets. Again, this does not come as a surprise given that for \TRPA the packing phase may even help to ``correct'' some choices made during the \rtprp phases. Regarding the \TRPA configurations, the number of extra packages required in the packing phase tends to be increasing in $\delta$. Thus, the transportation costs of \TPA and all \TRPA configurations get closer to each other after this phase. Interestingly, and somewhat surprisingly, we can see that the average transportation cost in Large-TS for the configuration with $\delta=0.85$ goes down by 4\% after the packing takes place. Again, this is consistent with the behavior observed in the plots of Figure~\ref{fig:obj_perf_round_pack} and, in particular, the improvement of the performance of the configuration with $\delta=0.85$ noted for Figure~\ref{fig:obj_large_pack}. Recall that Table~\ref{tab:packages} only accounts for the transportation costs and, therefore, although $\delta=0.85$ delivers the best such costs in Large-TS, this does not imply that it delivers the best costs overall. Indeed, Figure~\ref{fig:obj_large_pack} suggests that $\delta=0.95$ and $\delta=0.9$ are better when the penalization for unmet unvariable demand is also accounted for.

Summing up, the results presented in this section suggest that, when dealing with large problems that cannot be solved to optimality, \TRPA may be better than \TPA. Regarding the choice of parameter $\delta$, the above discussion does not offer a clear answer and specific analyses should be conducted for each potential application in order to ``calibrate'' this parameter.

\subsection{Results for Extra-large-TS}\label{sec:huge}

To conclude the numeric analysis we present a final comparison of the scaling of \tp and \rtp problems. The idea is to compare the percentage of instances in which they are able to obtain a feasible solution with two different time limits: 5 minutes and 30 minutes. Given that parameter $\delta$ does not change the difficulty of problem \rtp, we take $\delta=1$ in all executions. The analysis is developed independently for each of the 10~groups of instances of Extra-large-TS described in Section~\ref{sec:TS}. The size of the instances in each group has already been presented in Table~\ref{tab:huge_cases}.

The results, summarized in Table~\ref{tab:huge_cases_solved_instances}, show that \rtp can effectively obtain solutions in instances that have around twice the number of variables of those for which \tp can do the same. These results are consistent with the findings in the previous section, which already showed that \TRPA's behavior with respect to \TPA improved as problem instances became lager.

\begin{table}[!htbp]
    \centering
{\small
    \begin{tabular}{|c|c|cccccccccc|}
		\hline
    Time & Solution & \multirow{2}{*}{G1} & \multirow{2}{*}{G2} & \multirow{2}{*}{G3} & \multirow{2}{*}{G4} & \multirow{2}{*}{G5} & \multirow{2}{*}{G6} & \multirow{2}{*}{G7} & \multirow{2}{*}{G8} & \multirow{2}{*}{G9} & \multirow{2}{*}{G10} \\
		limit & method & & & & & & & & & & \\
    \hline
    \multirow{2}{*}{5 min.}  & \tp  & 100 & 100 & 100 & 100 & 100 &  64  &  8  & 0 & 0  & 0 \\                  
		& \rtp   & 100 & 100 & 100 & 100 & 100 & 100  & 96  & 4 & 4 & 0 \\               
		\hline                                    
    \multirow{2}{*}{30 min.} & \tp  & 100 & 100 & 100 & 100 & 100 & 100 & 100 & 100 & 80  & 64 \\
		& \rtp     & 100 & 100 & 100 & 100 & 100 & 100 & 100 & 100 & 100 & 100 \\
    \hline
		{\scriptsize Probl. Size} & {\scriptsize$\#$ vars.} & {\scriptsize$\approx 0.5\cdot 10^6$} & {\scriptsize $\times 1.9$} & {\scriptsize $\times 3.3$} & {\scriptsize $\times 5.2$} & {\scriptsize $\times 7.8$} & {\scriptsize $\times 11.1$} & {\scriptsize $\times 15.1$} & {\scriptsize $\times 20.2$} & {\scriptsize $\times 26.1$} & {\scriptsize $\times 33.2$} \\
		\hline
    \end{tabular}
		}
    \caption{Percentage of solved instances by \tp and \rtp for increasingly larger instances.}
    \label{tab:huge_cases_solved_instances}
\end{table}

\section{Conclusions}
\label{sec:conclusions}

We have presented a stock redistribution model for two-echelon supply networks. The model incorporates the packing component, which is very relevant for our partner in this project. Further, we believe it should also be relevant for most small and medium size retail companies and, in general, for any company which does not have its own fleet of transport vehicles. One of the main features of the model is its versatility, which, as we have discussed in Section~\ref{sec:policies}, allows to assess and compare different redistribution schemes.

Given the complex nature of the underlying mixed integer linear programming problem, designing efficient solution schemes is important in order to be able to apply it to real problems. With this goal in mind, we have developed an approximate algorithmic procedure, \TRPA, and we showed in Section~\ref{sec:TRP_num} that, for large instances, it outperforms the direct solution scheme, \TPA.

A natural step for future research would be to compare the performance, in large-scale problems, of \TRPA with decomposition techniques such as Branch and Price \citep{Vanderbeck1996,Barnhart1998,Vanderbeck2011} and Benders \citep{Benders1962,Rahmaniani2017}. Yet, it might be even more promising to study how the \TRPA scheme might be embedded in the above decomposition techniques so that one can improve the running times on the master problem and/or the subproblems.

\bibliography{Optimizacion.bib,StockManagement.bib}

\begin{thebibliography}{50}
\newcommand{\enquote}[1]{``#1''}
\expandafter\ifx\csname natexlab\endcsname\relax\def\natexlab#1{#1}\fi

\bibitem[\protect\citeauthoryear{Abouee{-}Mehrizi, Berman, and
  Sharma}{Abouee{-}Mehrizi et~al.}{2015}]{AboueeMehrizi2015}
\textsc{Abouee{-}Mehrizi, H., O.~Berman, and S.~Sharma} (2015):
  \enquote{Optimal Joint Replenishment and Transshipment Policies in a
  Multi-Period Inventory System with Lost Sales,} \emph{Operations Research},
  63, 342--350.

\bibitem[\protect\citeauthoryear{Agrawal, Chao, and Seshadri}{Agrawal
  et~al.}{2004}]{Agrawal2004}
\textsc{Agrawal, V., X.~Chao, and S.~Seshadri} (2004): \enquote{Dynamic
  balancing of inventory in supply chains,} \emph{European Journal of
  Operational Research}, 159, 296 -- 317, supply Chain Management: Theory and
  Applications.

\bibitem[\protect\citeauthoryear{Ahuja, Magnanti, and Orlin}{Ahuja
  et~al.}{1993}]{Ahuja1993}
\textsc{Ahuja, R.~K., T.~L. Magnanti, and J.~B. Orlin} (1993): \emph{Network
  Flows: Theory, Algorithms, and Applications}, Prentice hall.

\bibitem[\protect\citeauthoryear{Allen}{Allen}{1958}]{Allen1958}
\textsc{Allen, S.~C.} (1958): \enquote{Redistribution of total stock over
  several user locations,} \emph{Naval Research Logistics Quarterly}, 5,
  337--345.

\bibitem[\protect\citeauthoryear{Allen}{Allen}{1961}]{Allen1961}
\textsc{Allen, S.~G.} (1961): \enquote{A Redistribution Model with Set-Up
  Charge,} \emph{Management Science}, 8, 99--108.

\bibitem[\protect\citeauthoryear{Allen}{Allen}{1962}]{Allen1962}
---\hspace{-.1pt}---\hspace{-.1pt}--- (1962): \enquote{Computation for the
  Redistribution Model with Set-Up Charge,} \emph{Management Science}, 8,
  482--489.

\bibitem[\protect\citeauthoryear{Bacharach}{Bacharach}{1966}]{Bacharach1966}
\textsc{Bacharach, M.} (1966): \enquote{Matrix Rounding Problems,}
  \emph{Management Science}, 9, 732--742.

\bibitem[\protect\citeauthoryear{Banerjee, Burton, and Banerjee}{Banerjee
  et~al.}{2003}]{Banerjee2003}
\textsc{Banerjee, A., J.~Burton, and S.~Banerjee} (2003): \enquote{A simulation
  study of lateral shipments in single supplier, multiple buyers supply chain
  networks,} \emph{International Journal of Production Economics}, 81-82, 103
  -- 114, proceedings of the Eleventh International Symposium on Inventories.

\bibitem[\protect\citeauthoryear{Barnhart, Johnson, Nemhauser, Savelsbergh, and
  Vance}{Barnhart et~al.}{1998}]{Barnhart1998}
\textsc{Barnhart, C., E.~L. Johnson, G.~L. Nemhauser, M.~W.~P. Savelsbergh, and
  P.~H. Vance} (1998): \enquote{Branch-and-price: Column generation for solving
  huge integer programs,} \emph{Operations Research}, 46, 316--329.

\bibitem[\protect\citeauthoryear{Benders}{Benders}{1962}]{Benders1962}
\textsc{Benders, J.~F.} (1962): \enquote{Partitioning procedures for solving
  mixed-variables programming problems,} \emph{Numerische Mathematik}, 4,
  238--252.

\bibitem[\protect\citeauthoryear{Bendoly}{Bendoly}{2004}]{Bendoly2004}
\textsc{Bendoly, E.} (2004): \enquote{Integrated inventory pooling for firms
  servicing both on-line and store demand,} \emph{Computers \& Operations
  Research}, 31, 1465 -- 1480.

\bibitem[\protect\citeauthoryear{Bernhard and Vygen}{Bernhard and
  Vygen}{2006}]{Bernhard2006}
\textsc{Bernhard, K. and J.~Vygen} (2006): \emph{Combinatorial Optimization.
  Theory and Algorithms}, Algorithms and Combinatorics 21. Springer, Berlin,
  Heidelberg.

\bibitem[\protect\citeauthoryear{Burton and Banerjee}{Burton and
  Banerjee}{2005}]{Burton2005}
\textsc{Burton, J. and A.~Banerjee} (2005): \enquote{Cost-parametric analysis
  of lateral transshipment policies in two-echelon supply chains,}
  \emph{International Journal of Production Economics}, 93-94, 169 -- 178,
  proceedings of the Twelfth International Symposium on Inventories.

\bibitem[\protect\citeauthoryear{Caro and Gallien}{Caro and
  Gallien}{2010}]{Caro2010a}
\textsc{Caro, F. and J.~Gallien} (2010): \enquote{Inventory Management of a
  Fast-Fashion Retail Network,} \emph{Operations Research}, 58, 257--273.

\bibitem[\protect\citeauthoryear{Caro, Gallien, D{\'{\i}}az, Garc{\'{\i}}a,
  Corredoira, Montes, Ramos, and Correa}{Caro et~al.}{2010}]{Caro2010b}
\textsc{Caro, F., J.~Gallien, M.~D{\'{\i}}az, J.~Garc{\'{\i}}a, J.~M.
  Corredoira, M.~Montes, J.~A. Ramos, and J.~Correa} (2010): \enquote{Zara Uses
  Operations Research to Reengineer Its Global Distribution Process,}
  \emph{Interfaces}, 40, 71--84.

\bibitem[\protect\citeauthoryear{Christensen, Khan, Pokutta, and
  Tetali}{Christensen et~al.}{2017}]{Christensen2017}
\textsc{Christensen, H.~I., A.~Khan, S.~Pokutta, and P.~Tetali} (2017):
  \enquote{Approximation and online algorithms for multidimensional bin
  packing: {A} survey,} \emph{Comput. Sci. Rev.}, 24, 63--79.

\bibitem[\protect\citeauthoryear{Coffman, Csirik, Galambos, Martello, and
  Vigo}{Coffman et~al.}{2013}]{Coffman2013}
\textsc{Coffman, E.~G., J.~Csirik, G.~Galambos, S.~Martello, and D.~Vigo}
  (2013): \emph{Handbook of Combinatorial Optimization}, Springer New York,
  chap. Bin packing approximation algorithms: Survey and classification,
  455--531.

\bibitem[\protect\citeauthoryear{Conejo, Castillo, Minguez, and
  Garcia-Bertrand}{Conejo et~al.}{2006}]{Conejo2006}
\textsc{Conejo, A.~J., E.~Castillo, R.~Minguez, and R.~Garcia-Bertrand} (2006):
  \emph{Decomposition Techniques in Mathematical Programming}, Springer.

\bibitem[\protect\citeauthoryear{Cox and Ernst}{Cox and Ernst}{1982}]{Cox1982}
\textsc{Cox, L. and L.~Ernst} (1982): \enquote{Controlled Rounding,}
  \emph{INFOR: Information Systems and Operational Research}, 20, 423--432.

\bibitem[\protect\citeauthoryear{Diks and {de Kok}}{Diks and {de
  Kok}}{1996}]{Diks1996}
\textsc{Diks, E.~B. and A.~G. {de Kok}} (1996): \enquote{Controlling a
  divergent 2-echelon network with transshipments using the consistent
  appropriate share rationing policy,} \emph{International Journal of
  Production Economics}, 45, 369 -- 379, proceedings of the Eighth
  International Symposium on Inventories.

\bibitem[\protect\citeauthoryear{Diks and {de Kok}}{Diks and {de
  Kok}}{1998}]{Diks1998}
---\hspace{-.1pt}---\hspace{-.1pt}--- (1998): \emph{Lecture Notes in Economics
  and Mathematical Systems: {A}dvances in Distribution Logistics}, Springer,
  Berlin, Heidelberg, vol. 460, chap. Transshipments in a divergent 2-echelon
  system, 423–448.

\bibitem[\protect\citeauthoryear{Doerr, Friedrich, Klein, and Osbild}{Doerr
  et~al.}{2006}]{Doerr2006}
\textsc{Doerr, B., T.~Friedrich, C.~Klein, and R.~Osbild} (2006):
  \enquote{Unbiased Matrix Rounding,} in \emph{Algorithm Theory -- SWAT 2006},
  ed. by L.~Arge and R.~Freivalds, Berlin, Heidelberg: Springer Berlin
  Heidelberg, 102--112.

\bibitem[\protect\citeauthoryear{Dolan and Mor\'e}{Dolan and
  Mor\'e}{2002}]{Dolan2002}
\textsc{Dolan, E.~D. and J.~J. Mor\'e} (2002): \enquote{Benchmarking
  optimization software with performance profiles,} \emph{Mathematical
  Programming}, 91, 201--213.

\bibitem[\protect\citeauthoryear{Firoozi, Babai, Klibi, and Ducq}{Firoozi
  et~al.}{2020}]{Firoozi2020}
\textsc{Firoozi, M., M.~Z. Babai, W.~Klibi, and Y.~Ducq} (2020):
  \enquote{Distribution planning for multi-echelon networks considering
  multiple sourcing and lateral transshipments,} \emph{International Journal of
  Production Research}, 58, 1968--1986.

\bibitem[\protect\citeauthoryear{Geoffrion}{Geoffrion}{1972}]{Geoffrion1972}
\textsc{Geoffrion, A.~M.} (1972): \enquote{Generalized {B}enders
  decomposition,} \emph{Journal of Optimization Theory and Application}, 10,
  237--260.

\bibitem[\protect\citeauthoryear{Gross}{Gross}{1963}]{Gross1963}
\textsc{Gross, D.} (1963): \emph{Multistage Inventory Models and Techniques},
  Stanford University Press, Stanford, CA, chap. Centralized inventory control
  in multilocation supply systems, 47--84.

\bibitem[\protect\citeauthoryear{{Gurobi Optimization}}{{Gurobi
  Optimization}}{2020}]{GurobiOptimization2020}
\textsc{{Gurobi Optimization}} (2020): \enquote{Gurobi Optimizer Reference
  Manual,} Available at: \url{http://www.gurobi.com}.

\bibitem[\protect\citeauthoryear{Herer, Tzur, and Y\"ucesan}{Herer
  et~al.}{2006}]{Herer2006}
\textsc{Herer, Y.~T., M.~Tzur, and E.~Y\"ucesan} (2006): \enquote{The
  multilocation transshipment problem,} \emph{IIE Transactions}, 38, 185--200.

\bibitem[\protect\citeauthoryear{Hoadley and Heyman}{Hoadley and
  Heyman}{1977}]{Hoadley1977}
\textsc{Hoadley, B. and D.~P. Heyman} (1977): \enquote{A two-echelon inventory
  model with purchases, dispositions, shipments, returns and transshipments,}
  \emph{Naval Research Logistics Quarterly}, 24, 1--19.

\bibitem[\protect\citeauthoryear{Hu, Watson, and Schneider}{Hu
  et~al.}{2005}]{Hu2005}
\textsc{Hu, J., E.~Watson, and H.~Schneider} (2005): \enquote{Approximate
  solutions for multi-location inventory systems with transshipments,}
  \emph{International Journal of Production Economics}, 97, 31--43.

\bibitem[\protect\citeauthoryear{J\"onsson and Silver}{J\"onsson and
  Silver}{1987}]{Jonsson1987}
\textsc{J\"onsson, H. and E.~A. Silver} (1987): \enquote{Analysis of a
  Two-Echelon Inventory Control System with Complete Redistribution,}
  \emph{Management Science}, 33, 215--227.

\bibitem[\protect\citeauthoryear{Karmarkar}{Karmarkar}{1981}]{Karmarkar1981}
\textsc{Karmarkar, U.~S.} (1981): \enquote{The Multiperiod Multilocation
  Inventory Problem,} \emph{Operations Research}, 29, 215--228.

\bibitem[\protect\citeauthoryear{Karmarkar and Patel}{Karmarkar and
  Patel}{1977}]{Karmarkar1977}
\textsc{Karmarkar, U.~S. and N.~R. Patel} (1977): \enquote{The one-period,
  N-location distribution problem,} \emph{Naval Research Logistics Quarterly},
  24, 559--575.

\bibitem[\protect\citeauthoryear{Krishnan and Rao}{Krishnan and
  Rao}{1965}]{Krishnan1965}
\textsc{Krishnan, K.~S. and V.~R.~K. Rao} (1965): \enquote{Inventory control in
  {N} warehouses,} \emph{Journal of Industrial Engineering}, 212–215.

\bibitem[\protect\citeauthoryear{Lee, Jung, and Jeon}{Lee
  et~al.}{2007}]{Lee2007}
\textsc{Lee, Y.~H., J.~W. Jung, and Y.~S. Jeon} (2007): \enquote{An effective
  lateral transshipment policy to improve service level in the supply chain,}
  \emph{International Journal of Production Economics}, 106, 115 -- 126,
  special section on contextualisation of supply chain networks.

\bibitem[\protect\citeauthoryear{Lodi, Martello, and Monaci}{Lodi
  et~al.}{2002}]{Lodi2002}
\textsc{Lodi, A., S.~Martello, and M.~Monaci} (2002): \enquote{Two-dimensional
  packing problems: A survey,} \emph{European Journal of Operational Research},
  141, 241 -- 252.

\bibitem[\protect\citeauthoryear{Noham and Tzur}{Noham and
  Tzur}{2014}]{Noham2014}
\textsc{Noham, R. and M.~Tzur} (2014): \enquote{The single and multi-item
  transshipment problem with fixed transshipment costs,} \emph{Naval Research
  Logistics}, 61, 637--664.

\bibitem[\protect\citeauthoryear{Paterson, Kiesm{\"{u}}ller, Teunter, and
  Glazebrook}{Paterson et~al.}{2011}]{Paterson2011}
\textsc{Paterson, C., G.~P. Kiesm{\"{u}}ller, R.~H. Teunter, and K.~D.
  Glazebrook} (2011): \enquote{Inventory models with lateral transshipments:
  {A} review,} \emph{European Journal of Operational Research}, 210, 125--136.

\bibitem[\protect\citeauthoryear{Paterson, Teunter, and Glazebrook}{Paterson
  et~al.}{2012}]{Paterson2012}
\textsc{Paterson, C., R.~H. Teunter, and K.~D. Glazebrook} (2012):
  \enquote{Enhanced lateral transshipments in a multi-location inventory
  system,} \emph{European Journal of Operational Research}, 221, 317--327.

\bibitem[\protect\citeauthoryear{Rahmaniani, Crainic, Gendreau, and
  Rei}{Rahmaniani et~al.}{2017}]{Rahmaniani2017}
\textsc{Rahmaniani, R., T.~G. Crainic, M.~Gendreau, and W.~Rei} (2017):
  \enquote{The {B}enders decomposition algorithm: A literature review,}
  \emph{European Journal of Operational Research}, 259, 801--817.

\bibitem[\protect\citeauthoryear{Salazar-Gonz\'alez}{Salazar-Gonz\'alez}{2006}]{Salazar2006}
\textsc{Salazar-Gonz\'alez, J.} (2006): \enquote{Controlled rounding and cell
  perturbation: statistical disclosure limitation methods for tabular data,}
  \emph{Mathematical Programming}, 105, 585--603.

\bibitem[\protect\citeauthoryear{Tagaras and Vlachos}{Tagaras and
  Vlachos}{2002}]{Tagaras2002}
\textsc{Tagaras, G. and D.~Vlachos} (2002): \enquote{Effectiveness of stock
  transshipment under various demand distributions and nonnegligible
  transshipment times,} \emph{Production and Operations Management}, 11,
  183--198.

\bibitem[\protect\citeauthoryear{Tang and Yan}{Tang and Yan}{2010}]{Tang2010}
\textsc{Tang, S.-L. and H.~Yan} (2010): \enquote{Pre-distribution vs.
  post-distribution for cross-docking with transshipments,} \emph{Omega}, 38,
  192 -- 202.

\bibitem[\protect\citeauthoryear{Van~Rossum and Drake}{Van~Rossum and
  Drake}{2009}]{python3}
\textsc{Van~Rossum, G. and F.~L. Drake} (2009): \emph{Python 3 Reference
  Manual}, Scotts Valley, CA: CreateSpace.

\bibitem[\protect\citeauthoryear{Vanderbeck}{Vanderbeck}{2000}]{Vanderbeck2000}
\textsc{Vanderbeck, F.} (2000): \enquote{On Dantzig-Wolfe Decomposition in
  Integer Programming and ways to Perform Branching in a Branch-and-Price
  Algorithm,} \emph{Operations Research}, 48, 111--128.

\bibitem[\protect\citeauthoryear{Vanderbeck}{Vanderbeck}{2011}]{Vanderbeck2011}
---\hspace{-.1pt}---\hspace{-.1pt}--- (2011): \enquote{Branching in
  branch-and-price: a generic scheme,} \emph{Mathematical Programming}, 130,
  249--294.

\bibitem[\protect\citeauthoryear{Vanderbeck and Wolsey}{Vanderbeck and
  Wolsey}{1996}]{Vanderbeck1996}
\textsc{Vanderbeck, F. and L.~A. Wolsey} (1996): \enquote{An exact algorithm
  for IP column generation,} \emph{Operations Research Letters}, 19, 151--159.

\bibitem[\protect\citeauthoryear{Wee and Dada}{Wee and Dada}{2005}]{Wee2005}
\textsc{Wee, K.~E. and M.~Dada} (2005): \enquote{Optimal Policies for
  Transshipping Inventory in a Retail Network,} \emph{Management Science}, 51,
  1519--1533.

\bibitem[\protect\citeauthoryear{Wen, Choi, and Chung}{Wen
  et~al.}{2019}]{Wen2019}
\textsc{Wen, X., T.-M. Choi, and S.-H. Chung} (2019): \enquote{Fashion retail
  supply chain management: A review of operational models,} \emph{International
  Journal of Production Economics}, 207, 34 -- 55.

\bibitem[\protect\citeauthoryear{Zhang}{Zhang}{2005}]{Zhang2005}
\textsc{Zhang, J.} (2005): \enquote{Transshipment and Its Impact on Supply
  Chain Members’ Performance,} \emph{Management Science}, 51, 1534--1539.

\end{thebibliography}
\bibliographystyle{ecta}

\newpage

\small
\appendix
\section*{Appendix}
\label{sec:appendix}
\addcontentsline{toc}{section}{Appendices}
\renewcommand{\thesubsection}{\Alph{subsection}}

\subsection{Instance-generator pseudocode}
\label{sec:appendix_generator}

In this section we present a brief description of the pseudocode used for the generation of the instances of the different test sets studied in the paper. In the discussion that follows, the name of the parameter in the ensuing pseudocode procedures appears in parentheses.

The instance simulator allows to control various parameters in order to generate instances with the desired features for the numerical analysis at hand. We distinguish two main types of parameters:

\begin{description}
\item[Parameters related to instance sizes:] Number of SKUs ($num\_refs:=\card{\sku}$), different types of packages ($num\_packs:=\card{\pack}$), number of outlets ($num\_outlets := \card{\outl}$), and total network stock ($total\_stock:=\sum_{i\in\fac s\in\sku}\sini_{is}$). Moreover, in all simulations there is only one warehouse ($num\_warehouses := \card{\ware} = 1$).

\item[Parameters related to model features:] Redistribution policy ($movement\_policy\in \{\text{CR}, \text{DR}, \text{GR}\}$) and a scaling factor to capture different balances of the shipping rates of packages involving the warehouses with respect to direct shipments between outlets ($ware\_pack\_cost\_factor$).

\item[Parameters related to the solution methods:] Penalization parameters $\alpha$ and $\varepsilon$, along with $\delta$ parameter for the \rtp model, are not part of the instance generation procedure. Instead, they are parameters that can vary for different runs of \TPA and \TRPA algorithms on the same baseline instance. In particular, $\varepsilon$ is taken to be $0.0001$ in all runs, and $\alpha$ and $\delta$ are properly discussed in the body of the paper.
\end{description}

In all simulations it is assumed that 40\% of the stock is at the warehouse ($warehouses\_prop = 0.4$). The weight of the SKUS is between 0 and 1 ($min\_weight = 0$, $max\_weight = 1$). The baseline shipping costs are between 10 and 100 ($min\_cost = 10$, $max\_cost = 100$), and they are then perturbed depending on the involved outlets and type of package ($m\_factor = 0.5$, $mp\_factor = 0.8$). Package capacity will vary between 2 and 10 ($min\_cap = 2$, $max\_cap = 10$). For each SKU, the proportions of fixed demand and variable demand with respect to the total stock will vary between $0.5$ and 1 ($min\_fix\_dem\_factor = 0.5$, $max\_fix\_dem\_factor = 1$) and between 0.25 y 0.5 ($min\_var\_dem\_factor = 0.25$, $max\_var\_dem\_factor = 0.5$), respectively. Moreover, priorities $\pri_{is}$ are assumed to be equal across outlets ($min\_priority = max\_priority = 1$). 

%Por último, para la reproducibilidad de las simulaciones se establece una semilla que permite replicar las instancias generadas.

Procedure~\ref{alg:instances_generator_algorithm} presents the pseudocode of the instance generator. First, it defines the sets on which the mathematical model builds upon: $\sku$, $\pack$, $\ware$, $\outl$, $\fac$, and $\mov$. Next, it defines the model parameters. Parameters associated with the weights of the different SKUs ($\wei$) and package capacities ($\capa$) are generated from a uniform distribution. Parameters for costs ($\cost$), initial stocks  ($\sini$), fixed demand ($\fixd$), variable demand ($\vard$), and priorities ($\pri$) are defined in Procedures~\ref{alg:cost}-\ref{alg:priority} since some care is needed, in particular, to ensure feasibility of the resulting instances.

\begin{algorithm}[!htbp]
\caption{Generation of the \emph{cost} parameters ($\cost$)}
\label{alg:cost}
\begin{algorithmic}[1]
%\algrestore{generator}
\Require $\pack$; $\mov$
\Require $\capa$
\Require $min\_cost$; $max\_cost$; $m\_factor$; $mp\_factor$; $ware\_pack\_cost\_factor$

\State $max\_cap \gets \max_{p\in\pack}(\capa_{p})$
\State $ini\_min\_cost \gets min\_cost + (max\_cost-min\_cost) \cdot m\_factor \cdot mp\_factor$
\For{$p$ in $\Set{P}$}
    \State $ini\_cost_p \gets ini\_min\_cost + (max\_cost-ini\_min\_cost) \cdot \capa_p/max\_cap$
\EndFor
\For{$(i, j)$ in $\Set{M}$}
    \State $m\_cost \gets Uniform(m\_factor, 1)$
    \For{$p$ in $\Set{P}$}
        \State $move\_pack\_factor \gets Uniform(mp\_factor, 1)$
        \State $cost \gets move\_pack\_factor \cdot m\_cost \cdot ini\_cost_p$
        \If{$i$ or $j$ in $\Set{W}$}
            \State $\Param{Cost}_{ijp} \gets cost \cdot ware\_pack\_cost\_factor$
        \Else
            \State $\Param{Cost}_{ijp} \gets cost$
        \EndIf
    \EndFor
\EndFor
\Return $\Param{Cost}$
\end{algorithmic}
\end{algorithm}

\begin{algorithm}[!htbp]
\caption{Generation of the \emph{initial stock} parameters ($\sini$)}
\label{alg:ini_stock}
\begin{algorithmic}[1]
\Require $\sku$; $\ware$; $\outl$; $\fac$
\Require $\sini$
\Require $total\_stock$; $warehouses\_prop$

\State $t\_stock\_ware \gets Ceil(total\_stock \cdot warehouses\_prop)$
\State $t\_stock\_outlets \gets Floor(total\_stock  (1 - warehouses\_prop))$
\For{$i$ in $\fac$}
    \For{$s$ in $\sku$}
        \State $f\_stock_{is} \gets Uniform(0,1)$
    \EndFor
\EndFor
\For{$i$ in $\fac$}
    \For{$s$ in $\sku$}
        \If{$i$ in $\ware$}
            \State $sum\_f\_stock\_ware_{is} \gets f\_stock_{is} / \sum_{i^{\prime} \in \ware, s^{\prime} \in \sku} f\_stock_{i^{\prime}s^{\prime}}$
        \Else
            \State $sum\_f\_stock\_outlets_{is} \gets f\_stock_{is} / \sum_{i^{\prime} \in \outl, s^{\prime} \in \sku} f\_stock_{i^{\prime}s^{\prime}}$
        \EndIf
    \EndFor
\EndFor

\LeftComment{Warehouses}
\State $stock\_ware \gets Partitioning(t\_stock\_ware, sum\_f\_stock\_ware)^{\ \dagger}$%\footnotemark

\For{$i$ in $\ware$}
    \For{$s$ in $\sku$}
        \State $\sini_{is} \gets stock\_ware_{is}$
    \EndFor
\EndFor

\LeftComment{Outlets}
\State $stock\_outlets \gets Partitioning(t\_stock\_outlets, sum\_f\_stock\_outlets)$

\For{$i$ in $\outl$}
    \For{$s$ in $\sku$}
        \State $\sini_{is} \gets stock\_outlets_{is}$
    \EndFor
\EndFor
\Return $\sini$

\smallskip
$^\dagger${\scriptsize The $Partitioning(a, b)$ function splits $a$ into smaller taking into account the weights given by $b$. It is similar to multiply $a$ for each weight defined in $b$ but then ensuring that the sum of the partitioning is equal to $a$.}
\end{algorithmic}
\end{algorithm}

%\footnotetext{The $Partitioning(a, b)$ function splits $a$ into smaller taking into account the weights given by $b$. It is similar to multiply $a$ for each weight defined in $b$ but then ensuring that the sum of the partitioning is equal to $a$.}

\begin{algorithm}[!htbp]
\caption{Generation of the \emph{fixed demand} parameters ($\fixd$)}
\label{alg:fix_dem}
\begin{algorithmic}[1]
\Require $\sku$; $\ware$; $\outl$; $\fac$
\Require $\sini$
\Require $min\_fix\_dem\_factor; max\_fix\_dem\_factor$

\LeftComment{We need to do it by references to ensure that fix demand is satisfied}

\For{$s$ in $\sku$}
    \For{$i$ in $\outl$}
        \State $f\_fix\_dem_{is} \gets Uniform(0,1)$
    \EndFor
    
    \For{$i$ in $\outl$}
        \State $sum\_f\_fix\_dem_{is} \gets f\_fix\_dem_{is} / \sum_{i^{\prime}\in\outl} f\_fix\_dem_{i^{\prime}s}$
    \EndFor

    \State $t\_stock \gets \sum_{i \in \fac} \sini_{is}$
    \State $min\_fix\_dem \gets t\_stock \cdot min\_fix\_dem\_factor$
    \State $max\_fix\_dem \gets t\_stock \cdot max\_fix\_dem\_factor$
    \State $rand\_stock \gets Round(Uniform(min\_fix\_dem, max\_fix\_dem))$
    \State $fix\_dem \gets Partitioning(rand\_stock, sum\_f\_fix\_dem)$

    \For{$i$ in $\fac$}
        \If{$i$ in $\ware$}
            \State $\fixd_{is} \gets 0$
        \Else
            \State $\fixd_{is} \gets fix\_dem_{is}$
        \EndIf
    \EndFor
\EndFor
\Return $\fixd$
\end{algorithmic}
\end{algorithm}

\begin{algorithm}[!htbp]
\caption{Generation of the \emph{variable demand} parameters ($\vard$)}
\label{alg:var_dem}
\begin{algorithmic}[1]
\Require $\sku$; $\ware$; $\outl$; $\fac$
\Require $min\_var\_dem\_factor; max\_var\_dem\_factor$
\For{$i$ in $\outl$}
    \For{$s$ in $\sku$}
        \State $f\_var\_dem_{is} \gets Uniform(0,1)$
    \EndFor
\EndFor

\For{$i$ in $\outl$}
    \For{$s$ in $\sku$}
        \State $sum\_f\_var\_dem_{is} \gets f\_var\_dem_{is} / \sum_{i^{\prime}\in\outl, s^{\prime} \in \sku} f\_var\_dem_{i^{\prime}s^{\prime}}$
    \EndFor
\EndFor

\State $min\_var\_dem \gets total\_stock \cdot min\_var\_dem\_factor$
\State $max\_var\_dem \gets total\_stock \cdot max\_var\_dem\_factor$
\State $rand\_stock \gets Round(Uniform(min\_var\_dem, max\_var\_dem))$
\State $var\_dem \gets Partitioning(rand\_stock, sum\_f\_var\_dem)$
\For{$s$ in $\sku$}
    \For{$i$ in $\fac$}
        \If{$i$ in $\ware$}
            \State $\vard_{is} \gets 0$
        \Else
            \State $\vard_{is} \gets var\_dem_{is}$
        \EndIf
    \EndFor
\EndFor
\Return $\vard$
\end{algorithmic}
\end{algorithm}

\begin{algorithm}[!htbp]
\caption{Generation of the \emph{priority} parameters ($\pri$)}
\label{alg:priority}
\begin{algorithmic}[1]
\Require $\sku$; $\ware$; $\outl$; $\fac$
\Require $min\_priority; max\_priority$
\For{$i$ in $\fac$}
    \For{$s$ in $\sku$}
        \If{$i$ in $\ware$}
            \State $\pri_{is} \gets 0$
        \Else
            \State $\pri_{is} \gets Uniform(min\_priority, max\_priority)$
        \EndIf
    \EndFor
\EndFor
\Return $\pri$
\end{algorithmic}
\end{algorithm}

\begin{algorithm}[!htbp]
\caption{Instance-generator algorithm}
\label{alg:instances_generator_algorithm}
\begin{algorithmic}[1]
\Require $num\_refs > 0$; $num\_packs > 0$; $num\_outlets > 0$; $total\_stock > 0$ 
\Require $movement\_policy \in \{\textrm{CR},\ \textrm{DR},\ \textrm{GR}\}$
\Require $ware\_pack\_cost\_factor \geq 0$
\State $num\_warehouses \gets 1$; $warehouses\_prop \gets 0.4$
\State $min\_weight \gets 0$; $max\_weight \gets 1$
\State $min\_cost \gets 10$; $max\_cost \gets 100$
\State $m\_factor \gets 0.5$; $mp\_factor \gets 0.8$
\State $min\_cap \gets 2$; $max\_cap \gets 10$
\State $min\_fix\_dem\_factor \gets 0.5$; $max\_fix\_dem\_factor \gets 1$
\State $min\_var\_dem\_factor \gets 0.25$; $max\_var\_dem\_factor \gets 0.5$
\State $min\_priority \gets 1$; $max\_priority \gets 1$

\AlgSection{Sets}
\State $\sku \gets \{1,\ldots,num\_refs\}$
\State $\pack \gets \{1,\ldots,num\_packs\}$
\State $\ware \gets \{0\}$
\State $\outl \gets \{1,\ldots,num\_outlets\}$
\State $\fac \gets \ware \cup \outl$
\If{$movement\_policy = \textrm{CR}$}
    \State $\Set{M} \gets \{(i, j)\in\fac\times\fac: i\neq j\ \mathrm{and}\ (i \in \ware\ \textrm{or}\ j \in \ware)\}$
\ElsIf{$movement\_policy = \textrm{DR}$}
    \State $\Set{M} \gets \{(i, j)\in\fac\times\fac: i\neq j\ \mathrm{and}\ j \in \outl\}$
\Else
    \State $\Set{M} \gets \{(i, j)\in\Set{S}\times\Set{S}: i\neq j\}$
\EndIf

\AlgSection{Parameters}
\For{$s$ in $\sku$}
    \State $\wei_i \gets Uniform(min\_weigth, max\_weight)$
\EndFor
\For{$p$ in $\pack$}
    \State $\capa_p \gets Uniform(min\_cap, max\_cap)$
\EndFor
\State $\cost \gets \mathrm{Procedure\ \ref{alg:cost}}$
\State $\sini \gets \mathrm{Procedure\ \ref{alg:ini_stock}}$
\State $\fixd \gets \mathrm{Procedure\ \ref{alg:fix_dem}}$
\State $\vard \gets \mathrm{Procedure\ \ref{alg:var_dem}}$
\State $\pri \gets \mathrm{Procedure\ \ref{alg:priority}}$
\end{algorithmic}
\end{algorithm}

\end{document}